\documentclass[11pt]{article}
\usepackage{amsmath}
\usepackage{graphicx}
\usepackage{epsfig}
\usepackage{amsfonts}
\usepackage{amssymb}
\usepackage{placeins}
\usepackage{cases}
\usepackage[margin=1in]{geometry}
\usepackage{authblk}
\usepackage[titletoc,toc,title]{appendix}
\usepackage{epstopdf}
\usepackage{changes}
\usepackage{xcolor}
\usepackage{bm}
\usepackage{enumerate}

\def\horizontaldistance{\kern2pt}
\def\verticaldistance{\kern 5pt}

\usepackage{accents}

\title{Reduced Order Modeling of Nonlinear Dynamical Systems Using Slow Manifolds}

\begin{document}
\author[1]{Dan Wilson \thanks{corresponding author:~dwilso81@utk.edu}}
\affil[1]{Department of Electrical Engineering and Computer Science, University of Tennessee, Knoxville, TN 37996, USA}
\maketitle

\begin{abstract}
Model order reduction in high-dimensional, nonlinear dynamical systems if often enabled through fast-slow timescale separation.  One such approach involves identifying a low-dimensional slow manifold to which the state rapidly converges and subsequently studying the behavior on the slow manifold.  This work investigates slow manifolds defined by the intersection of an unstable manifold of an unstable fixed point or periodic orbit and the stable manifold of a stable attractor.  When the decay rates of perturbations transverse to the unstable manifold are sufficiently large, the resulting slow manifold can be used for reduced modeling purposes by leveraging the isostable coordinate framework. Detailed examples are provided for two different highly nonlinear dynamical systems, the first being a coupled system of Hodgkin-Huxley neurons and the second being a  biophysically detailed model of circadian oscillations.  The resulting reduced order models are illustrated in  two different biologically motivated control objectives. 
\end{abstract}

\section{Introduction}

High dimensionality is a common factor that precludes mathematical analysis and control design in many nonnegligibly nonlinear dynamical systems necessitating the use of model order reduction algorithms as a preliminary step in the analysis.  Model order reduction of nonlinear dynamical systems is often enabled by exploiting an inherent timescale separation between the fast and slow dynamics.  For instance, center manifold theory \cite{wigg03} can be applied in certain situations (and is particularly useful for bifurcation analysis).  In a similar manner, inertial manifolds \cite{foia88}, \cite{foia882} can be used to define invariant manifolds for dynamical systems that attract solutions exponentially quickly.  Singular perturbation theory can be applied when separating fast variables from slow variables to understand how the dynamical behavior evolves on a slow manifold  \cite{feni79}, \cite{kape99} \cite{farj18}, \cite{guck09}.  The dimension of linear models that capture the dynamics near a stable fixed point can be reduced by first diagonalizing the system and truncating all rapidly decaying eigenmodes, as gauged by the real component of the associated eigenvalues.  This approach can be extended to nonlinear dynamical systems using spectral submanifolds \cite{hall16}, \cite{pons20}, \cite{cene22} which are often approximated through asymptotic expansion.

In many applications, isostable coordinates can be used to study dynamical behavior that is characterized by timescale separation between fast an slow dynamics.  Isostable coordinates can be formally defined as the principal Koopman eigenmodes \cite{mezi20}, \cite{kval21} associated with either a stable fixed point or a periodic orbit.  For instance, considering a general dynamical system 
\begin{equation} \label{maineqintro}
    \dot{x} = F(x),
\end{equation}
with $x\in \mathbb{R}^N$, letting $x_0$ be a fixed point with $\lambda_1,\dots,\lambda_N$ being eigenvalues associated with its linearization, for any solution $x(t)$ that satisfies \eqref{maineqintro}, the isostable coordinates $\psi_1,\dots,\psi_N$ evolve in time according to
\begin{equation}
    \frac{d \psi_j(x(t))}{dt} = \lambda_j \psi_j,
\end{equation}
 for $j = 1,\dots,N$  \cite{wils20ddred}, \cite{wils21dd}.  Ordering the isostable coordinates in terms of their decay rate so that $| {\rm Real}(\lambda_k)| \leq |{\rm Real}(\lambda_{k+1})|$, provided there is some separation between the real components of the eigenvalues, a small subset of the slowest decaying isostable coordinates can be used to define a reduced order coordinate system \cite{maur13}, \cite{wils16isos}, \cite{wils17isored} that can subsequently be used for mathematical analysis and control design  \cite{wils21input}, \cite{ahme23}, \cite{park24}. 

 Isostable coordinate-based approaches are similar to Koopman-based approaches such as dynamic mode decomposition (DMD) \cite{schm10}, \cite{will15}, \cite{kutz16} or strategies that seek to identify Koopman invariant subspaces \cite{brun16}, \cite{lusc18}, \cite{kais21} attempt to provide an approximate linear representation of a nonlinear dynamical system using a lifted coordinate basis.  However, Isostable-coordinate-based methods differ from these approaches in that they retain the essential nonlinear characteristics of the underlying system, providing a coordinate system that isolates the slow dynamics.  Indeed, using isostable coordinates it is straightforward to define a slow manifold:
\begin{equation} \label{slowdefintro}
    S_{\rm iso} = \{ x \in \mathbb{R}^N | \psi_k(x) = 0 \; {\rm for} \; k > \beta  \}.
\end{equation}
Provided there is a large enough gap between $\psi_\beta$ and $\psi_{\beta+1}$, solutions rapidly converge to the slow manifold so that the dynamics of the full, nonlinear system can be well approximated by considering the dynamics on $S_{\rm iso}$.  While \eqref{slowdefintro} provides a relatively simple definition of a slow manifold, numerical challenges associated with the separation between fast and slow timescales often make computation of $S_{\rm iso}$ challenging.  Previous work has focused on strategies that compute an asymptotic expansion of the slow manifold in a close neighborhood of a fixed point \cite{wils21dd}, \cite{maur16} or a periodic orbit \cite{wils17isored}, \cite{wils20highacc}.  These approaches, however, are typically insufficient for use with model order reduction when large magnitude inputs are required and the computational complexity of these approaches makes them difficult to implement in high dimensional systems.  Recent work \cite{wils25} investigated strategies for computation of trajectories in backward time along the slow manifold, but these methods were difficult to implement when computing the slow manifold far from the underlying stable attractor.  

\begin{figure}[htb]
\begin{center}
\includegraphics[height=3.0 in]{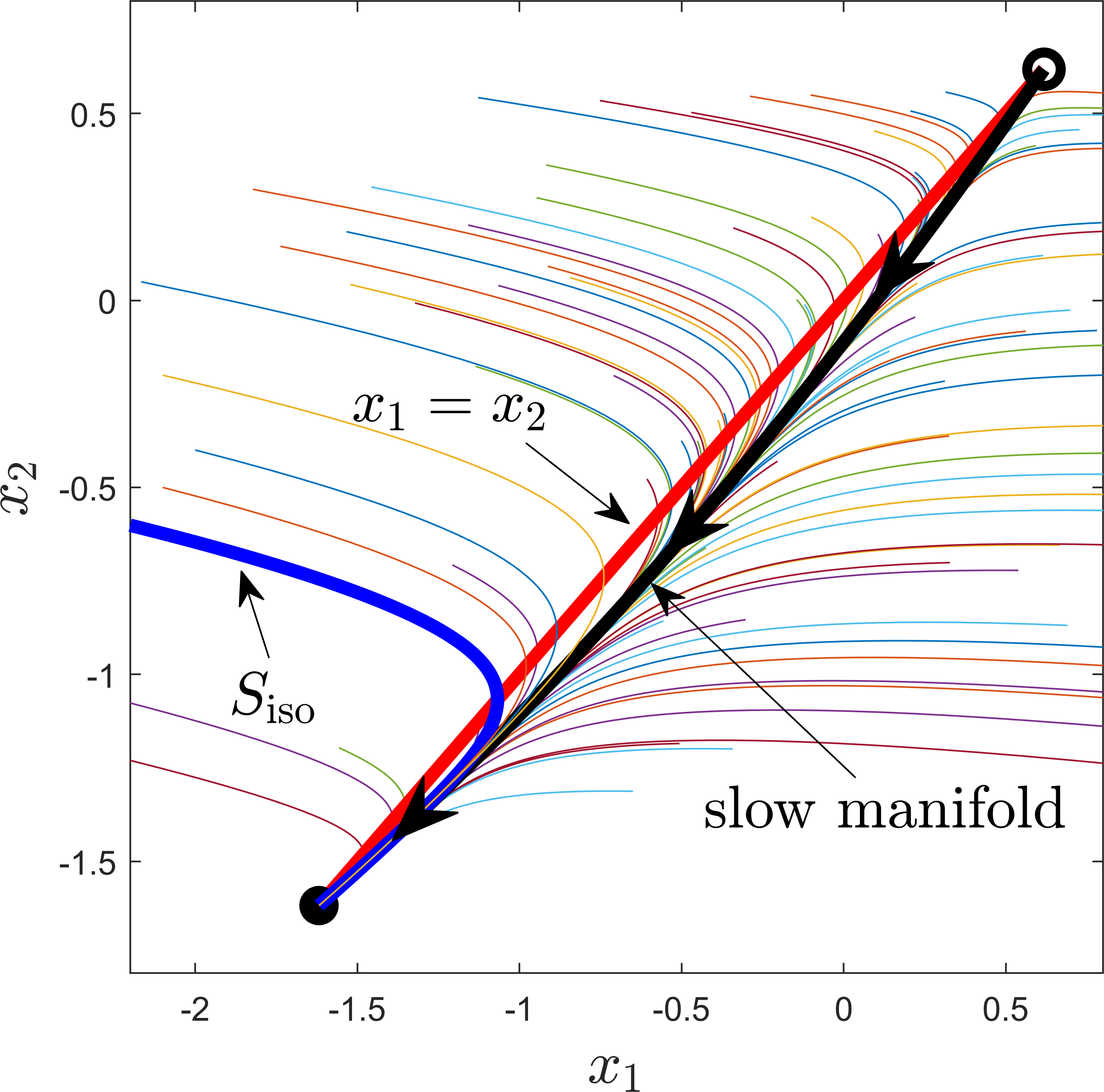}
\end{center}
\caption{Thin colored lines show 100 example trajectories evolving under the flow of \eqref{toymodel}.  The stable and unstable fixed point of this model are shown with solid and open circles, respectively.  States rapidly collapse to a slow manifold, given by the intersection of the stable manifold of the stable fixed point and the unstable manifold of the unstable fixed point.  The line $x_1 = x_2$ is shown for comparison, which approximates the slow manifold, but is not an invariant manifold.  The blue line defined according to \eqref{slowdefintro}  does not behave as a slow manifold in this example system.}
\label{simplemodelfigure}
\end{figure}

In contrast to previously developed computational methods for the identification of slow manifolds, this work explores the slow manifold in relation to the stable/unstable fixed points and periodic orbits of the underlying dynamical system \eqref{maineqintro}.  As a toy example, consider a 2-dimensional  model 
\begin{align} \label{toymodel}
\dot{x}_1 & = -\mu(x_1-x_2), \nonumber \\
\dot{x}_2 & = -1 + x_2^2 + x_1.
\end{align}  
Here, $\mu = 8$.  The model \eqref{toymodel} has an unstable fixed point at $x_1 = x_2 = 0.5 + \sqrt 5/2$ and a stable fixed point at $x_1 = x_2 = 0.5 - \sqrt 5/2$.   For large values of $\mu$, trajectories evolving under the flow of the vector field rapidly collapse to points near the line $x_1 = x_2$ (red line in Figure \ref{simplemodelfigure}).  The linearized stable fixed point has eigenvalues $\lambda_1 = -1.92$ and $\lambda_2 = -9.32$.   The manifold $S_{\rm iso}$  defined by \eqref{slowdefintro} is approximated to 12th order with an asymptotic expansion using techniques described in \cite{wils21dd} and is shown in blue.  Near the fixed point, this manifold does indeed behave as a slow manifold, but eventually veers to the left.  It is immediately apparent that trajectories are rapidly attracted to the black curve, which is the intersection of the unstable manifold of the unstable fixed point and the stable manifold of the stable fixed point.  

With the example from Figure \ref{simplemodelfigure} in mind, this manuscript investigates slow manifolds that are defined by the intersection of an unstable manifold of an unstable fixed point or periodic orbit and the stable manifold of a stable attractor and the use of such slow manifolds for reduced order modeling purposes.  The organization of this paper is as follows:~Section \ref{backsec} provides necessary background information on isostable coordinates of both stable fixed points and periodic orbits.  Section \ref{proofsection} considers strategies for characterizing the decay rate of solutions to the unstable manifold of both fixed points and periodic orbits.  Section \ref{ressec} provides two numerical illustrations, identifying a slow manifold in a coupled system of Hodgkin-Huxley neurons and in a biophysically detailed model for circadian oscillations.  In each example, the resulting reduced order model is used to formulate and solve two biologically motivated control objectives.  Section \ref{concsec} provides concluding remarks.  



\FloatBarrier

\section{Background} \label{backsec}

\subsection{Isostable Coordinates of Fixed Points}

Consider a general dynamical system
\begin{equation} \label{maineq}
    \dot{x} = F(x) + u(t),
\end{equation} 
where $x \in \mathbb{R}^N$, $F$ gives the dynamics, and $u$ is an external input.  For the moment, taking $u = 0$, suppose \eqref{maineq} has  with a stable fixed point $x_0$ for which $F(x_0) = 0$.  Local linearization yields
\begin{equation} \label{lineq}
   \Delta \dot{x} = J \Delta x + O(||\Delta x||^2),
\end{equation}
where $\Delta x = x - x_0$ and $J$ is the Jacobian evaluated at $x_0$.  While \eqref{lineq} is only valid in a close neighborhood of the fixed point, the spectrum of $F$ can be used to define isostable coordinates which correspond to level sets of principal Koopman eigenfunctions \cite{mezi19}, \cite{maur13}, \cite{mezi13} and can ultimately be used characterize the behavior of solutions in the basin of attraction of the fixed point.  To do so, let $w_k$, $v_k$ and $\lambda_k$ be left eigenvectors, right eigenvectors, and eigenvalues of $J$, respectively, ordered so that $|{\rm Real}(\lambda_k)| \leq |{\rm Real}(\lambda_{k+1})|$.  Additionally, suppose that the left and right eigenvectors are scaled so that $w_k^T v_k = 1$ where $^T$ denotes the transpose.  For the slowest decaying eigenvalue, $\lambda_1$, an associated principal isostable coordinate $\psi_1(x)$ can be defined according to
\begin{equation} \label{isodef}
    \psi_1(x) = \lim_{t \rightarrow \infty} (w_1^T (\phi(t,x) - x_0) \exp (- \lambda_1 t)),
\end{equation}
where $\phi(t,x)$ denotes the flow of \eqref{maineq} when $u(t) = 0$.  Additional isostable coordinates can be defined as level sets of principal Koopman eigenfunctions \cite{kval21}.  When $u(t) = 0$, under the flow of \eqref{maineq} in the basin of attraction of the fixed point, isostable coordinates evolve in time according to 
\begin{equation} \label{psidyn}
    \dot{\psi}_j = \lambda_j \psi_j,
\end{equation}
for $j = 1,\dots,N$.  To a linear approximation, the isostable coordinates are equivalent to using an eigenbasis associated with the local linearization, i.e.,~provided $\psi_k = O(\epsilon)$ where $0 < \epsilon \ll 1$
\begin{equation} \label{stateapprox}
    x = x_0 + \sum_{k = 1}^\beta \psi_k v_k + O(\epsilon^2).
\end{equation}
As distinct from local linearization, however, as shown in \cite{wils21dd}, it is possible extend the approximation \eqref{stateapprox} to higher orders of accuracy.

\subsection{Phase and Isostable Coordinates of Periodic Orbits}
Phase and isostable coordinates can also be defined for systems of the form \eqref{maineq} that have a stable $T$-periodic orbit $x^\gamma(t)$ when $u = 0$.  The phase coordinates capture the timing of oscillations while isostable coordinates give a sense of the decay of perturbations transverse to the limit cycle.  A phase $\theta \in [0,2\pi)$ can be assigned for all $x \in x^\gamma$, scaled so that $d \theta/dt = \omega$ for solutions evolving under the flow of \eqref{maineq}.  Isochrons \cite{winf01}, \cite{guck75} can be used to define phase in the basin of attraction of $x^\gamma$ as follows:~for some $a(0)\in x^\gamma$, letting $\theta(a(0)) = \theta_1$, the $\theta_1$ level set (i.e.,~isochron) is defined as the set of all $b(0)$ such that 
\begin{equation}
    \lim_{t \rightarrow \infty} || a(t) - b(t)|| = 0,
\end{equation}
where $|| \cdot||$ is some vector norm.  Isostable coordinates for periodic orbits can also be defined leveraging Floquet theory \cite{jord07}.  Local linearization near the periodic orbit yields
\begin{equation} \label{perlin}
    \Delta \dot{x} = J(t) \Delta x + O(||\Delta x||^2),
\end{equation}
where $\Delta x = x-x^\gamma(t)$ and $J(t)$ is the Jacobian evaluated at $x^\gamma(t)$.  Letting $\Phi$ be the fundamental matrix of the linear time-varying system \eqref{perlin} for which $\Delta x (T) = \Phi \Delta x(0)$.   For simplicity of exposition, suppose that $\Phi$ is diagonalizable with eigenvalues $\lambda_0, \lambda_1, \dots, \lambda_{N-1}$.  Floquet theory allows solutions of \eqref{perlin} to be written as
\begin{equation} \label{floqsoln}
    \Delta x(t) = \sum_{j = 0}^{N-1} c_j \exp(\kappa_j t) p_j(t) + O(||\Delta x||^2),
\end{equation}
where $\kappa_j = \log(\lambda_j)/T$ are Floquet exponents, $p_j(t)$ are $T$-periodic Floquet eigenfunctions, and $c_j$ are constants that depend on initial conditions.  The Floquet exponents are often sorted so that $\kappa_0 = 0$ (associated with the Floquet multiplier $\lambda_0 = 1$).  The remaining Floquet exponents can be sorted according to ${\rm Real}(\kappa_j) \geq {\rm Real}(\kappa_{j+1})$.  Provided $x^\gamma$ is stable, the slowest decaying isostable coordinate can be defined explicitly according to (cf.,~\cite{wils16isos})
\begin{equation} \label{iso1def}
    \psi_1 = \lim_{k \rightarrow \infty} \big[ q_1^T(t) (x(t) - x^\gamma(t)) \exp(- \kappa_1 t) \big],
\end{equation}
where $q_1^T(t) p_1(t) = 1$.  Much like the definition \eqref{isodef}, which encodes for the slowest decaying eigenmode of a system with a stable fixed point, \eqref{iso1def} encodes for the slowest decaying Floquet eigenmode as solutions decay to the stable periodic orbit.  Additional isostable coordinates $\psi_2,\dots,\psi_{N-1}$ can be defined as level sets of principal Koopman eigenfunctions \cite{kval21}, \cite{mezi20} of the stable periodic orbit.  When $u(t) = 0$, under the flow of \eqref{maineq}, in the basin of attraction of the periodic orbit, phase and isostable coordinates evolve in time according to
\begin{align} \label{psidynper}
    \dot{\theta} &= \omega, \nonumber \\
    \dot{\psi}_j &= \kappa_j \psi_j,
\end{align}
for $j = 1,\dots,N-1$.  Notice that isostable coordinates have the same unperturbed decay for systems with periodic orbits \eqref{psidyn} and for systems with fixed points \eqref{psidynper}.

\subsection{Dynamics of Isostable Coordinates for Fixed Points} \label{fpdynsec}

Isostable coordinates have previously been used in a wide variety of reduced order modeling applications,  \cite{maur16}, \cite{soot17}, \cite{wils21dd}, \cite{wils21input}, \cite{ahme23}, \cite{wils20highacc}, \cite{wils21adapt}.  For isostable coordinates defined for both fixed point and periodic orbits, the evolution of the isostable coordinates along solutions $x(t)$ of Equation \eqref{maineq} is
\begin{align} \label{isodyn}
    \frac{d \psi_j}{dt} &= \frac{\partial \psi_j}{\partial x}^T \frac{dx}{dt} \nonumber \\
    &= I_j^T (F(x) + u(t)) \nonumber \\
    &= \lambda_j \psi_j + I_j^T u(t),
\end{align}
for $j = 1,\dots,N$.  Above, $I_j \equiv \frac{\partial \psi_j}{\partial x}$ evaluated at $x(t)$ and the simplification in the third line results from the fact that $\dot{\psi}_j = \lambda_j \psi_j$ when $u = 0$.  Along unperturbed trajectories of \eqref{maineq},  $I_j(t)$ evolves according to \cite{wils20highacc}, 
\begin{equation} \label{isoeq}
\dot{I}_j = -(J^T - \lambda_j {\rm Id}) I_j,
\end{equation}
for $k = 1,\dots,N$, where $J$ is the Jacobian evaluated at $x(t)$ and ${\rm Id}$ is an appropriately sized identity matrix.  Noting that when $x(t) = x_0$,  $\dot{I}_j = 0$ so that Equation \eqref{isoeq} simplifies to $0 = (J^T - \lambda_j {\rm Id}) I_j$ with solution $I_j = w_jk$.  In situations where there is a large gap between the decay rates of the isostable coordinates, i.e.,~if ${\rm Real}(\lambda_\beta)$ is much larger than ${\rm Real}(\lambda_{\beta + 1})$ the faster decaying isostable coordinates can be assumed to be zero to arrive at a reduced order model.  Using the isostable reduced model \eqref{isodyn}, the state can be approximated as a function of the nonzero isostable coordinates.

\subsection{Dynamics of Phase and Isostable Coordinates for Periodic Orbits}
The dynamics of phase and isostable coordinates of periodic orbits are almost identical to those of isostable coordinates of fixed points, with the main difference being the inclusion of phase coordinates to encode for the timing of oscillations.  Along solutions $x(t)$ of Equation \eqref{maineq}, isostable coordinates evolve in time according to 
\begin{align} \label{isophasedyn1}
\frac{d \psi_j}{dt} &= \frac{\partial \psi_j}{\partial x}^T \frac{dx}{dt} \nonumber \\
&= I_j^T(F(x) + u(t)) \nonumber \\
&= \kappa_j \psi_j + I_j^T u(t).
\end{align}
Likewise, phase coordinates evolve in time according to
\begin{align} \label{isophasedyn2}
\frac{d \theta}{dt} &= \frac{\partial \theta}{\partial x}^T \frac{dx}{dt} \nonumber \\
&= Z^T(F(x) + u(t)) \nonumber \\
&= \omega + Z^T u(t).
\end{align}
Above, $I_j \equiv \frac{\partial \psi_j}{\partial x}$ and $Z \equiv \frac{\partial \theta}{\partial x}$, both evaluated at $x(t)$.  The simplification in the third line results from the fact that $\dot{\psi}_j = \kappa_j \psi_j$ and $\dot{\theta} = \omega$ when $u = 0$.  Along unperturbed  trajectories of \eqref{maineq}, $Z(t)$ and $I_j(t)$ evolve in time according to \cite{wils21dd}
\begin{align}  \label{zidotper}
    \dot{Z} &= -J^T Z, \nonumber \\
    \dot{I}_j &= -(J^T  - \kappa_j {\rm Id}) I_j,
\end{align}
for $j = 1,\dots,{N-1}$ where $J$ is the Jacobian evaluated at $x(t)$.  Above, the evolution equation for the isostable coordinate gradients are identical to those from \eqref{isoeq}.   As in the case with stable fixed points, if ${\rm Real}(\lambda_\beta)$ is much larger than ${\rm Real}(\lambda_{\beta + 1})$, the faster decaying isostable coordinates can be assumed to be zero to arrive at a reduced order model from wihch the state can be approximated as a function of the nonzero isostable coordinates.

\section{Slow Manifolds Connecting Unstable Equilibria and Stable Attractors} \label{proofsection}
For systems with either a stable fixed point or a stable periodic orbit, model order reduction is possible when $|{\rm Real}(\lambda_{\beta})|  -  |{\rm Real}(\lambda_{\beta+1})|   $ (or $|{\rm Real}(\kappa_{\beta})| - |{\rm Real}(\kappa_{\beta+1})|$) is large for some $\beta$.  In these instances, the fast isostable coordinates $\psi_{\beta+1}, \dots, \psi_N$ decay rapidly relative to the slow isostable coordinates so that can be well-approximated by zero and truncated, ultimately yielding a reduced order model.  With this in mind, using the isostable coordinate framework, a slow manifold, $S_{\rm iso}$ can be defined according to \eqref{slowdefintro} \cite{wils25}.  Provided there is a large gap between the real components of eigenvalues $\lambda_\beta$ and $\lambda_{\beta+1}$, solutions of \eqref{maineq} will converge to $S_{\rm iso}$ rapidly.  The fundamental challenge in implementing this strategy lies in the numerical computation of $S_{|rm iso
}$.  A naive approach would be to consider some initial condition near the fixed point for which $\psi_k(x) = 0$ for $k > \beta$ and integrate backwards in time along the slow manifold.  However, strongly attracting directions in forward time become strongly repelling directions in backward time making this approach numerically infeasible.  Reference \cite{wils25} suggested two approaches for computation of $S_{\rm iso}$, however, these were only able to provide an approximation of the slow manifold.  

Alternatively, one can consider a general dynamical system with an unstable fixed point $x_0$ that coexists with either a stable fixed point or periodic orbit $x_s$.  Let $W^u_{x_0}$ be the unstable manifold of an unstable fixed point and $W^s_{x_s}$ be the stable manifold of $x_s$.  As shown below, under appropriate technical conditions, 
\begin{equation} \label{slowstate1}
    S =  W_{x_0}^u \cap W^s_{x_s},
\end{equation}  
defines a slow manifold to which nearby states converge exponentially fast.  With Equation \eqref{slowstate1} in mind, one can approximate $S$ near $x_0$ through local linearization and simply compute the rest of the manifold along the slow manifold using forward time integration.   This slow manifold can subsequently be used to understand the dynamics of \eqref{maineq} in a reduced order setting.  
Likewise, for a system \eqref{maineq} with an unstable periodic orbit $x^\gamma$ and a stable attractor $x_s$,
\begin{equation} \label{slowstate2}
    S = W^u_{x^\gamma} \cap W_{x_s}^s,
\end{equation}
where $W_{x_s}^s$ is the stable manifold of $x_s$ and $W^u_{x^\gamma}$ is the unstable manifold of $x^\gamma$ also defines a slow manifold to which nearby states converge exponentially fast.  Once again, from Equation \eqref{slowstate2}, the slow manifold can be computed  by obtaining an approximation of $W^u_{x^\gamma}$ near the unstable periodic orbit and integrating forward in time.


\subsection{Globally Linearized Coordinates on the Unstable Manifold of an Unstable Fixed Point}
Consider the case for a system of the form \eqref{maineq} with an unstable fixed point $x_0$. For the moment, suppose $u(t) = 0$.  Local linearization can be used to capture the behavior when $x(t) - x_0 = O(\epsilon)$ where $0 < \epsilon \ll 1$
\begin{equation}
     x(t) - x_0 = \sum_{j= 0}^N  \alpha_k(t) v_k(t) + O(\epsilon^2).
\end{equation}
where $v_k$ is an eigenvector of $A = \frac{\partial F}{\partial x}$ evaluated $x_0$ and $\alpha_k$ is the coordinate in the basis of eigenvectors.  Provided $F$ is twice differentiable, Hartman-Groban theorem can be used to show that the nonlinear system \eqref{maineq} is at least locally conjugate to a linear system.  Work in \cite{lan13} extended this theorem to the basin of attraction of the unstable manifold of the fixed point.  Leveraging these results, letting $W^u_{x_0}$ be the unstable manifold of the fixed point, there exists $h(x):W_{x_0}^u \rightarrow \mathbb{C}^\beta$ such that $y = h(x)$ is a $C^1$ diffeomorphism with $\frac{\partial h}{\partial x}|_{x = 0} = {\rm Id}$ that satisfies 
\begin{equation} \label{linfp}
    \dot{y} = \Lambda y.
\end{equation}  
Here, $\Lambda$ is a diagonal matrix containing the unstable eigenvalues $\lambda_1,\dots,\lambda_\beta$ of $A$.  Note that the theorem proved in \cite{lan13} states $\dot{y} = Ay$, but the form stated above can be obtained through a subsequent diagonalization.  As a function of the coordinates ${y}$, one can write
\begin{align} \label{ucoords}
    x-x_0 &=   {f}({y}),
\end{align}
where $ {f}({y}) \equiv  {h}^{-1}({y})$.  Considering two infinitesimally close trajectories, $x_1,x_2 \in W^u_{x_0}$, $\Delta x = x_2 - x_1$ can be written to leading order as
\begin{align}
\Delta x &= \frac{\partial {f}}  {\partial {{y}}}  \Delta {y}, \nonumber \\
&= g_j({y}) \Delta {y}_j,
\end{align}
where $\Delta {y} = {h}(x_2)-{h}(x_1)$, $g_j({y})$ is the $j^{\rm th}$ column of $\frac{\partial {f}}  {\partial {{y}}}$, and $\Delta {y}_j$ is the $j^{\rm th}$ entry of $\Delta {y}$.  Noting that when ${y} = O(\epsilon)$, where $0 < \epsilon \ll 1$,
\begin{equation}
   \Delta x = \sum_{j = 1}^\beta v_j \Delta {y}_j + O(\epsilon).
\end{equation}
As such,
\begin{equation}
    g_j({y}) \approx v_j
\end{equation}
when ${y}$ is small.  One can approximate $g_j(y(t))$ along a given trajectory ${x}(t) \in W^u_{x_0}$ when $x(0)$ starts close to a fixed point according to 
\begin{equation} \label{findg}  
    g_j(y(t)) = \frac{x_1(t) -x(t)} {\epsilon \exp(\lambda_j t)},
\end{equation}
where $x_1(0) = {x}(0) + \epsilon v_j$.

\subsection{Globally Linearized Coordinates on the Unstable Manifold of an Unstable Periodic Orbit}
One can also consider a slow manifold for a system of the form \eqref{maineq} according to \eqref{slowstate2}, i.e.,~as the intersection the unstable manifold of an unstable periodic orbit and the stable manifold of a stable attractor.  Once again, for the moment let $u(t) = 0$.  As described in \cite{lan13} suppose that \eqref{maineq} has an unstable, $T$-periodic orbit $x^\gamma$ and let $W^u_{x^\gamma}$ be the $\beta$-dimensional unstable manifold of this periodic orbit.   By first rescaling time so that the period of the orbit is $2 \pi$, it is possible to define a new coordinate system in a neighborhood of the periodic orbit
\begin{align}
\dot{\theta} &= 1, \label{seq} \\
\dot{z} &= m(z,\theta),  \label{zeq}
\end{align}
with $\theta \in \mathbb{S}^1$ and $z \in \mathbb{R}^{\beta-1}$.  Here $\theta$ is a time-like variable in the longitudinal direction and $z$ captures the behavior transverse to the periodic orbit.   For $x \in x^\gamma$, $z = 0$.  As shown in \cite{lan13}, the $2\pi$-periodic system \eqref{zeq} is topologically conjugate to a system of the form 
\begin{equation} \label{jthetaeq}
    \dot{y} = J(\theta) y,
\end{equation} 
through a $C^1$ map $y = b(z,\theta)$; here $J(\theta) = \frac{\partial g}{\partial z}|_{z = 0}$.  As a consequence, the original system \eqref{maineq} is related to the linearized system comprised of \eqref{seq} and \eqref{jthetaeq} through the diffeomorphism $x \rightarrow (\theta(x),z(x)) \rightarrow (\theta(x),b(z(x),\theta(x)))$.  As discussed in \cite{lan13}, this diffeomorphism extends to the entire unstable manifold of the periodic orbit.

Following a similar argument as the case for stable fixed points.  As discussed above, \eqref{maineq} can be linearized by the diffeomorphism $(\theta,y) = (h_1(x),h_2(x))$ with inverse 
\begin{equation} \label{inverseeqn}
    x = H(y,\theta).
\end{equation} 
Here, $\theta\in [0, 2 \pi)$ and $y \in \mathbb{R}^{\beta-1}$.   In this transformed coordinate system, $\dot{\theta} = 1$ and $\dot{y} = J(\theta) y$.  Since $\theta$ is $2\pi$-periodic, Floquet theory \cite{jord07} can be used to write solutions in the transformed coordinate system as
\begin{align}
\theta(t) &= \theta(x(0)) + t \nonumber, \\
y(t) &= \sum_{j = 1}^{\beta-1} c_j(x(0)) \exp(\kappa_j t) p_j(\theta),
\end{align}
where $\kappa_j$ are Floquet exponents, $p_j(\theta)$ are $2 \pi$-periodic Floquet eigenfunctions, and $c_j$ are constants that depend on initial conditions.  Transforming back to the original coordinates, recalling that $x \in x^\gamma$ when $y = 0$, one can write
\begin{align} \label{perexp}
    x - x^\gamma(\theta) &= \frac{\partial H}{\partial y} y(t) + O(||y^2||) \nonumber \\
    &=  \sum_{j = 1}^{\beta-1} g_j(\theta) c_j(x(0))\exp(\kappa_j t) + O(||y^2||),
\end{align}
where $\frac{\partial H}{\partial y}$ is evaluated at $y = 0$,  $m_j(\theta) = \frac{\partial H}{\partial y} p_j(\theta)$.  Note that because the periodic orbit is unstable, $\rm{Real}(\kappa_j) >0$ for all $j$.  Near the periodic orbit, the terms $g_j(\theta)$ can be approximated with a local linearization  to yield
\begin{equation} \label{perlin}
    \Delta \dot{x} = J(t) \Delta x + O(||\Delta x||^2),
\end{equation}
where $\Delta x = x - x^\gamma(\theta)$, $J(t)$ is the Jacobian evaluated at $x^\gamma(t)$, and $\theta = \omega t$.  Assuming that the Fundamental matrix is diagonalizable, Floquet theory allows one to write solutions of \eqref{perlin} on the unstable manifold as
\begin{equation} \label{floqsoln}
    x - x^\gamma(\theta) = \sum_{j = 1}^{\beta-1}   p_j(\theta)  c_j(x(0)) \exp(\kappa_j t)  + O(||\Delta x||^2),
\end{equation}
where $\kappa_j$ are the unstable Floquet exponents, $p_j(\theta)$ are the corresponding $2\pi$-periodic Floquet eigenfunctions, and $c_j$ are constants that depend on initial conditions.  Comparing \eqref{floqsoln} and \eqref{perexp}, $g_j(\theta)$ and $p_j(\theta)$ are identical to leading order when $y$ and $\Delta x$ are small.  

Considering two infinitesimally close trajectories, $x_1,x_2 \in W_{x^\gamma}^u$, using \eqref{inverseeqn}, to leading order $\Delta x = x_2-x_1$ can be written as
\begin{align}
    \Delta x &= \frac{\partial H}{\partial y} \Delta y + \frac{\partial H}{\partial \theta} \Delta \theta, \nonumber \\
    &= \sum_{j = 1}^{\beta-1} g_j(y,\theta) \Delta y_j + g_\beta(y,\theta) \Delta \theta,
\end{align}
where $\Delta y = h_2(x_2)-h_2(x_1)$, $\Delta \theta = h_1(x_2)-h_2(x_1)$, $g_1,\dots,g_{\beta-1}$ is given by the $j^{\rm th}$ column of $\frac{\partial H}{\partial y}$, $g_{\beta} = \frac{\partial H}{\partial \theta}$, and $\Delta y _j$ is the $j^{\rm th}$ entry of $\Delta y$.  As stated above, $g_j(y,\theta)$ are well approximated by the Floquet eigenfunctions of the linearization when $y$ is small.  Likewise, from \eqref{perexp} one finds that to leading order, $g_\beta = \frac{\partial x^\gamma}{\partial \theta}$ when $y$ is small.  One can compute $g_j(y(t),\theta(t))$ for $j = 1,\dots,\beta-1$ along a given trajectory $x(t) \in W^u_{x^\gamma}$ when $x(0)$ starts close to the periodic orbit according to
\begin{equation} \label{gper1}
    g_j(y(t),\theta(t)) = \frac{x_1(t)-x(t)}{\epsilon \exp(\kappa_j t)},
\end{equation}
where $x_1(0) = x(0) + \epsilon g_j(\theta)$. Likewise, noting that $\Delta \theta$ does not change over time between the trajectories $x_2(t)$ and $x(t)$, one finds
\begin{equation}  \label{gper2}
    g_\beta(y(t),\theta(t)) = \frac{x_2(t)-x(t)}{\epsilon},
\end{equation}
where $x_2(0) = g_\beta(\theta)$.

\subsection{Local Decay of Trajectories to the Unstable Manifold}
The decay of trajectories to the unstable manifold will be considered here using the linearized coordinate system for states evolving on the unstable manifold of either a fixed point or periodic orbit as discussed in the previous sections.  To this end, let $x(t) = x_0 + {f}({y}(t))$ (resp.,~$x(t) = H(y(t),\theta(t))$) be a trajectory evolving on the unstable manifold of a fixed point (resp.,~periodic orbit).  Consider a perturbation $\Delta x$ transverse the unstable manifold, i.e.,~in the complement of ${\rm span}\{g_1,\dots,g_\beta\}$.  This perturbation can be written as
\begin{align}
\Delta x(t)&=    N(t) N(t)^\dagger \Delta x(t) + G(t) G(t)^\dagger \Delta x(t),
\end{align}
where $G(t) = \begin{bmatrix} g_1(t) & \dots &  g_\beta(t)  \end{bmatrix}$, $N$ is the orthogonal complement of $G$, and $^\dagger$ denotes the pseudoinverse. Above, the contribution $N(t) N(t)^\dagger \Delta x(t)$ represents the component that is transverse to the unstable manifold; as such, the decay to the slow manifold will be considered by focusing on the decay of this term.  Assuming that $||\Delta x|| = O(\epsilon)$ where $0<\epsilon \ll 1$, to leading order, 
\begin{equation} \label{deltaxeq}
    \Delta \dot{x} = J(t) \Delta x,
\end{equation}
where $J(t)$ is the Jacobian evaluated at $x(t)$.  Noting that \eqref{deltaxeq} is a linear time-varying equation, solutions of $\Delta x(t)$ follow
\begin{align}
    \Delta x(t) &= \Phi(t,0) \Delta x(0) \nonumber \\
    &=  N(t) N(t)^\dagger  \Phi(t,0)\Delta x(0)  +    G(t) G(t)^\dagger \Phi(t,0) \Delta x(0) \nonumber \\
    &= \Delta x_N(t) + \Delta x_G(t),
\end{align}
where $\Phi(t,0)$ is the state transition matrix associated with \eqref{deltaxeq} and $\Delta x_N(t)$ and $\Delta x_G(t)$ are defined appropriately.  While we are interested in the evolution of $||\Delta x_N(t)||_2$, this is numerically challenging because $||\Delta x_N(t)||_2$ generally shrinks at a far smaller rate than $||\Delta x_G(t)||_2$.  As a workaround to this issue, instead consider 
\begin{equation} \label{xhatdot}
\Delta \dot{\hat{x}} = J(t) \Delta \hat{x}  - \rho G(t) G(t)^\dagger   \Delta \hat{x},
\end{equation}
where $\rho$ is a positive constant.  First, notice that 
\begin{align}  \label{firstnn}
    \frac{d}{dt}  \big(N(t) N(t)^\dagger \Delta \hat{x} \big) &= \frac{d}{dt} \big(  N(t) N(t)^\dagger  \big) \Delta \hat{x} + N(t) N(t)^\dagger \Delta \dot{\hat{x}} \nonumber \\
    &=  \frac{d}{dt} \big(  N(t) N(t)^\dagger  \big) \Delta \hat{x} + N(t) N(t)^\dagger  \big(J(t) \Delta \hat{x}  - \rho G(t) G(t)^\dagger   \Delta \hat{x}) \big) \nonumber \\
    &= \frac{d}{dt} \big(  N(t) N(t)^\dagger  \big) \Delta \hat{x} + N(t) N(t)^\dagger  J(t) \Delta \hat{x},
\end{align}
where the second line is obtained by substituting \eqref{xhatdot} and the third line results from the fact that $N(t) N(t)^\dagger G(t) G(t)^\dagger$ is identical to the zero matrix.  Comparing with
\begin{align} \label{secondnn}
     \frac{d}{dt}  \big(N(t) N(t)^\dagger \Delta {x} \big) &= \frac{d}{dt} \big(  N(t) N(t)^\dagger  \big) \Delta {x} + N(t) N(t)^\dagger \Delta \dot{{x}} \nonumber \\
     &= \frac{d}{dt} \big(  N(t) N(t)^\dagger  \big) \Delta {x} + N(t) N(t)^\dagger  J(t) \Delta {x}, 
\end{align}
one finds that solutions of  \eqref{firstnn} are identical to those of \eqref{secondnn}.   Next considering the evolution of terms in the direction of the unstable manifold:
\begin{align}
\frac{d}{dt} \big(  G(t) G(t)^\dagger  \Delta \hat{x}  \big) &=  \frac{d}{dt} \big(  G(t) G(t)^\dagger  \big) \Delta \hat{x} + G(t) G(t)^\dagger \Delta \dot{\hat{x}}  \nonumber \\
&=  -\rho G(t) G(t)^\dagger \Delta \hat{x} + \bigg(  G(t) G(t)^\dagger J(t)  + \frac{d}{dt} \big( G(t) G(t)^\dagger  \big)  \bigg) \Delta \hat{x},
\end{align}
where the second line is obtained by substituting \eqref{xhatdot} and rearranging.  Noting that $G(t)G(t)^\dagger$ and $J(t)$ are predetermined by the trajectory $x(t)$, the contribution from $G(t) G(t)^\dagger  \Delta \hat{x} $ can be made arbitrarily small be choosing $\rho$ large enough.  Combining the above insights,  one finds that by choosing $\rho$ to be large enough \eqref{xhatdot} isolates the contribution of $\Delta x_N(t)$ from $\Delta x(t)$.  Finally, in order to mitigate issues related to the fact that $\Delta x_N(t)$ generally approaches zero exponentially fast, one can instead consider the solution of 
\begin{equation} \label{xbardot}
    \Delta \dot{{\bar x}} = J(t) \Delta \bar{x} - \rho G(t)G(t)^\dagger \Delta \bar{x} - \nu(t) \Delta \bar{x},
\end{equation}
where $\nu(t) \in \mathbb{R}$ is chosen so that  $\frac{d}{dt} ||\bar{x}||_2 = 0$, i.e.,~so that the norm of the solution $\bar{x}(t)$ does not change in time.  One can verify that solutions of \eqref{xhatdot} and \eqref{xbardot} are related by
\begin{equation}
    \Delta \hat{x} (t) = \Delta \bar{x}(t)   \left(   \int_0^t \nu(\tau) d\tau) \right).
\end{equation}
Here, $\eta(t)$ is taken to be the local decay rate, which can be obtained directly considering the constraint
\begin{align}
   0 &= \frac{d}{dt} || \Delta \bar{x} ||_2^2 \nonumber \\
   &= \frac{d}{dt} (\Delta \bar{x}^T)   \Delta \hat{x} + \Delta \bar{x}^T \frac{d}{dt}(\Delta \hat{x}) \nonumber \\
   &= \Delta \bar{x}^T \big( J(t)^T + J(t) + \rho ( G(t)G(t)^\dagger +  (G(t)G(t)^\dagger)^T )  \big) \Delta \bar{x}  - 2 v(t) \Delta \bar{x}^T \Delta \bar{x},
\end{align}
which can be rearranged to yield
\begin{equation} \label{nuequation}
    \nu(t) = \frac{\Delta \bar{x}^T \big( J(t)^T + J(t) + \rho ( G(t)G(t)^\dagger +  (G(t)G(t)^\dagger)^T )  \big) \Delta \bar{x} }{2 v(t) \Delta \bar{x}^T \Delta \bar{x}  }.
\end{equation}

To summarize, for a given trajectory $x(t)$ evolving on the $\beta$-dimensional unstable manifold, $g_1(t),\dots,g_\beta(t)$ can be approximated along trajectories starting close to an unstable fixed point using \eqref{findg} or an unstable periodic orbit using \eqref{gper1} and \eqref{gper2}.  The decay in directions transverse to the unstable manifold can be determined by considering the set $\{n_1(t),\dots,n_{N-\beta}(t)  \}$ with initial conditions $\{n_1(0),\dots,n_{N-\beta}(0)  \}$ chosen to be in the orthogonal complement of the span of $\{g_1(0),\dots,g_\beta(0)\}$.  Local decay rates, $\nu_k(t)$ associated with each $n_k(t)$ can be computed solving Equation \eqref{nuequation} with initial condition $n_k(0)$.  This process can be repeated with an appropriate sampling of trajectories on the unstable manifold to characterize decay rate of perturbations to the unstable manifold.

\section{Reduced Order Modeling in Numerical Examples} \label{ressec}

Once a slow manifold of the form \eqref{slowstate1} or \eqref{slowstate2} has been identified, it is possible to use a $\beta$-dimensional isostable coordinate-based reduced order model to approximate the full $N$-dimensional model dynamics \eqref{maineq}.  Using the results from Section \ref{proofsection} it is generally possible to identify the slow manifold by finding an approximation of the unstable manifold near the unstable fixed point (resp.,~periodic orbit) and propagating forward in time until it converges to the stable attractor.  Once this slow manifold has been found, a reduced order model of the form \eqref{isodyn} can be used to characterize the behavior on the slow manifold when the associated stable attractor is a fixed point.  When the stable attractor is a periodic orbit, a similar reduction can be used that uses both phase and isostable coordinates in a model of the form \eqref{isophasedyn1} and \eqref{isophasedyn2}.  The following section provides three examples that illustrate the utility of this approach.

\subsection{Toy Model Illustration}
\begin{figure}[htb]
\begin{center}
\includegraphics[height=2.7 in]{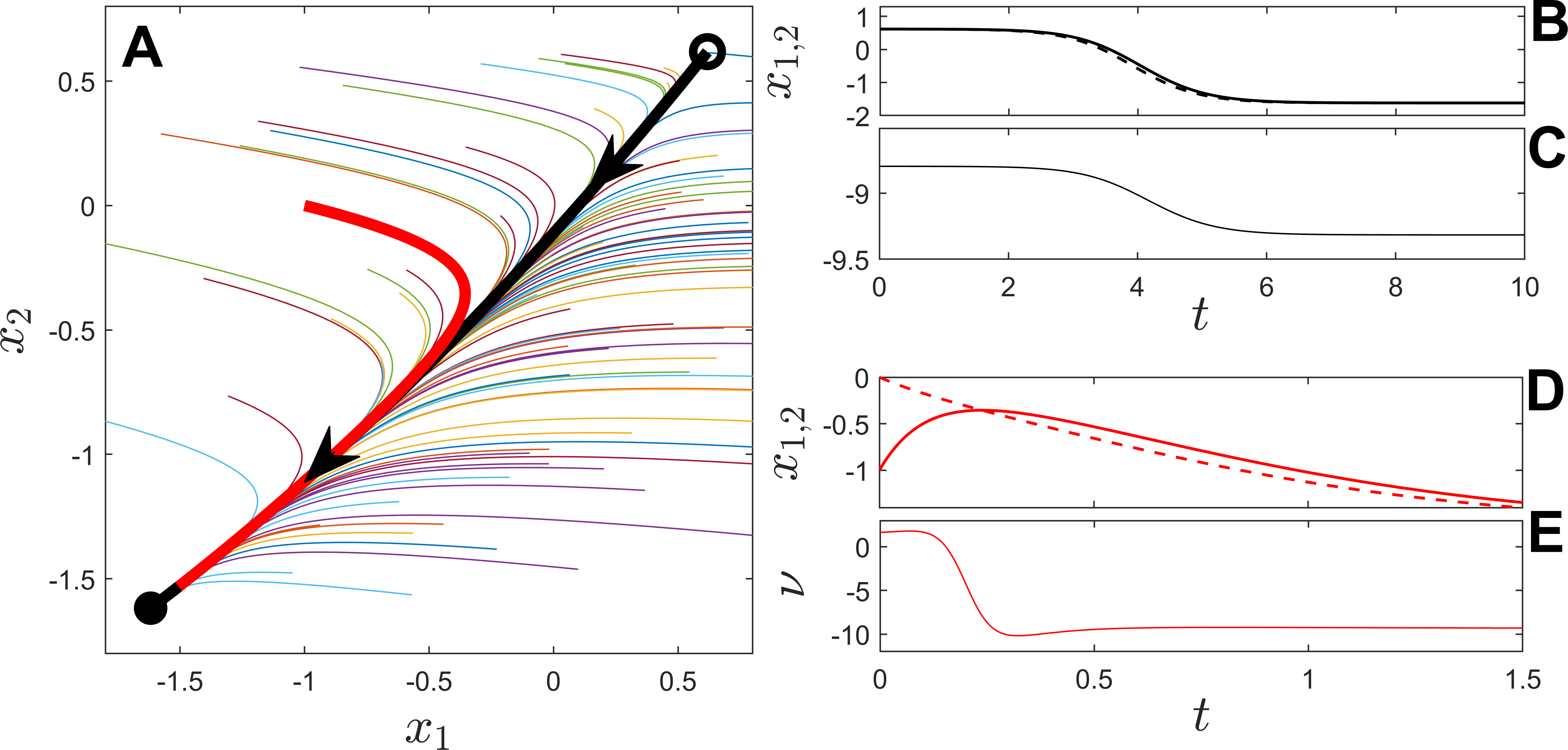}
\end{center}
\caption{For the toy model \eqref{toymodel}, Panel A shows the intersection of the unstable manifold of the unstable fixed point and the stable manifold of the stable fixed point in black.  Local trajectories converge rapidly to this slow manifold.  Panel B shows $x_1(t)$ and $x_2(t)$ for reference and Panel C shows $\nu(t)$ evaluated along this trajectory, confirming the rapid decay in transverse directions.  The analysis is repeated for the red curve, with panel D showing $x_1(t)$ and $x_2(t)$ along this curve for reference.  Panel E shows the corresponding value of $\nu(t)$ indicating that initially, nearby trajectories do not converge to the red curve so that it cannot be treated as a slow manifold.}
\label{simplemodelanalysis}
\end{figure}

Here, the behavior of the toy model \eqref{toymodel} introduced earlier is considered in terms of a slow manifold reduction.  As shown in panel A of Figure \ref{simplemodelanalysis}, trajectories rapidly converge to the intersection of the unstable manifold of the unstable fixed point and the stable manifold of the stable fixed point (black line).  In panel B, both $x_1$ and $x_2$ are shown along this curve as the solid and dashed lines, respectively with the decay rate of solutions transverse to the slow manifold $\nu(t)$ computed according to \eqref{nuequation} an shown in panel C.  Indeed, solutions near this slow manifold converge to it rapidly as evidenced by the negative value of $\nu$ all along the trajectory.  For the comparison, this analysis is repeated for the red trajectory in panel A which is initially far from the slow manifold before rapidly converging to it.  For this trajectory, $g_1 = \begin{bmatrix} \dot{x_1} & \dot{x_2}\end{bmatrix}^T$ represents the direction of travel of the solutions along this trajectory.  The decay rate of solutions transverse to the flow, $\nu(t)$, for the red trajectory is also computed according to \eqref{nuequation} and shown in Panel E.  Panel D shows $x_1(t)$ and $x_2(t)$ along the red curve for reference.  Initially along this trajectory, because $v(t)>0$, perturbations transverse to the red curve diverge before ultimately converging as the trajectory approaches the true slow manifold. 

\begin{figure}[htb]
\begin{center}
\includegraphics[height=2.4 in]{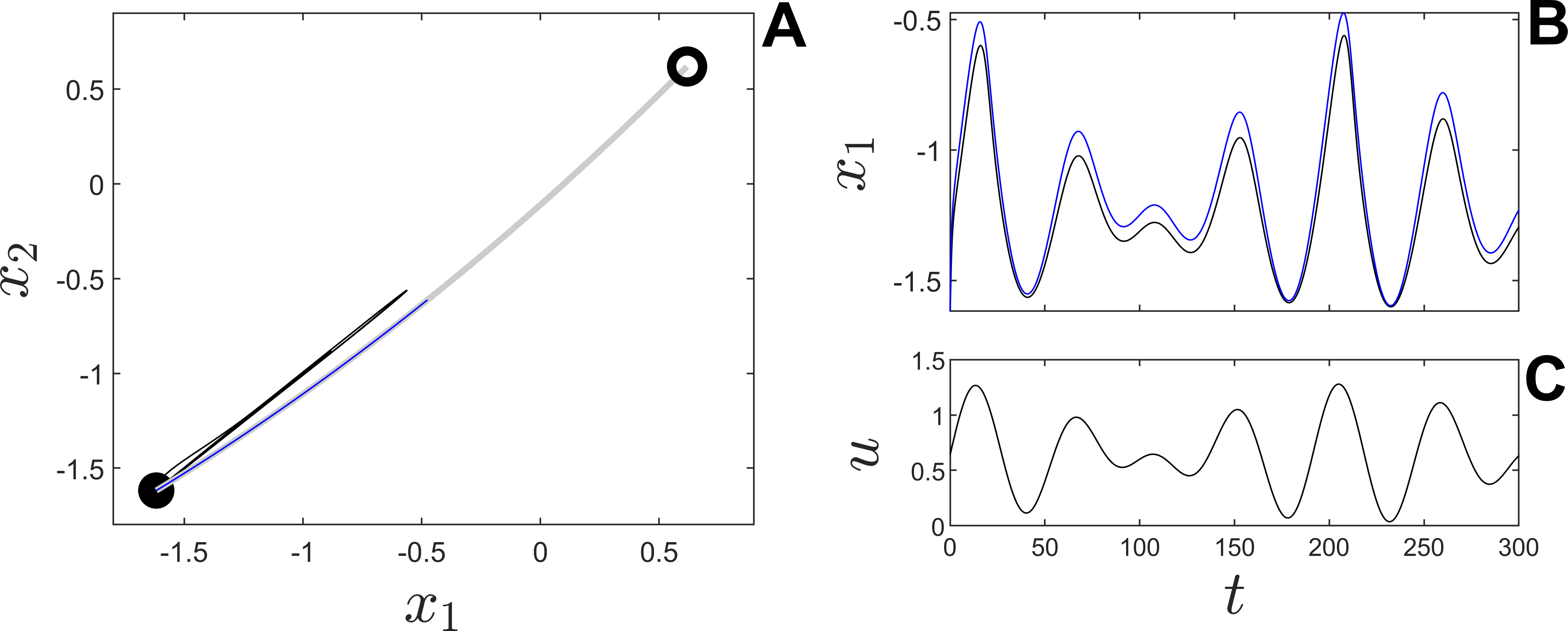}
\end{center}
\caption{The full model \eqref{toymodel} and the isostable coordinate-based reduced order represenation are both simulated in response to the input $u(t)$ in panel C.  Panel A shows the slow manifold in grey along with the stable and unstable fixed points shown as closed and open circles, respectively.  The black and blue lines show the full and reduced order model simulations, respectively.  Panel B shows $x_1$ plotted as a function of time for these two simulations.  As discussed in the text, the agreement between the two models improves as $\mu$ decreases.  }
\label{reducedsimple}
\end{figure}

The linearized stable fixed point has two eigenvalues $\lambda_1 = -1.92$ and $\lambda_2 = 9.32$.  A reduced order model of the form \eqref{isodyn} is obtained for the  model \eqref{toymodel} considering slow decaying isostable coordinate $\psi_1$ associated with $\lambda_1$ and truncating the contribution from fast decaying isostable coordinate.   $I_1$ is computed using Equation \eqref{isoeq} in backwards time along the slow manifold.  Figure \ref{reducedsimple} compares the evolution of the state from the isostable reduced and full order models in response to the input $u(t) = 0.64 + 0.32\sin(0.1 t) + 0.32 \sin(0.13 t)$.  Panel A shows the resulting trajectory of the reduced order model in the $x_1$-$x_2$ plane.  Panel B shows the evolution of $x_1$ over time for each model and panel $C$ shows the applied input for reference.  Letting $x_1(t)$ and $x_2(t)$ be the solution of \eqref{toymodel} and $\tilde{x}_1$ and $\tilde{x}_2$ be the solution obtained from the isostable reduced model, the average error over the course of this simulation is defined by
\begin{equation}
  E =   \frac{1}{300} \int_0^{300} ( (x_1-\tilde{x}_1)^2 + (x_2-\tilde{x}_2)^2 )^{1/2}
\end{equation}
and is equal to 0.0691.  This test is repeated for different values of $\mu$.  Solutions decay faster to the slow manifold as $\mu$ decreases, ultimately decreasing the error between the full and reduced order models.   For instance, $E = 0.0363$ when $\mu = -15$ and $E = 0.0238$ when $\mu = -25$.

\FloatBarrier

\subsection{Maintaining Quiescence in a Model of Coupled Hodgkin-Huxley Neurons} \label{hhsec}
Consider a model of four coupled Hodgkin-Huxley model neurons:
\begin{align} \label{hheqs}
     C_m  \dot{V}_j &= -I_{\rm ion}(V_j,n_j,m_j,h_j)  - g_C(V_j - \bar{V}) + B_j + \alpha_j u(t),
\end{align}
for $j = 1,\dots,4$.   Here $V_j$ is the transmembrane voltage of neuron $j$ in mV,  $I_{\rm ion}(t)$ is comprised of a set of ionic currents mediated by gating variables $n_j$, $m_j$, and $h_j$,  $C_m = 1 \mu{\rm F}/{\rm cm}^2$ is the membrane capacitance,  $B_j = -4.318 + 0.75j$ $\mu {\rm A}/{\rm cm}^2$ is the baseline current of each neuron, and time is in milliseconds.  Coupling is electrotonic \cite{john95} with $g_C = 1.5$ and $\bar{V} = \frac{1}{4} \sum_{j = 1}^4 V_j$.  Here $u(t)$ is an input in the form of a transmembrane current applied to only the first two neurons so that $\alpha_j = 1$ for $j = 1,2$ and equals zero for $j = 3,4$.  Note that synaptic coupling could also be used instead of electrotonic coupling if desired.  A full set of model equations is given in Appendix \ref{hhapx}.  With a total of 4 coupled neurons, each with 4 state variables, this model has a total of 16 state variables.

\begin{figure}[htb]
\begin{center}
\includegraphics[height=2.0 in]{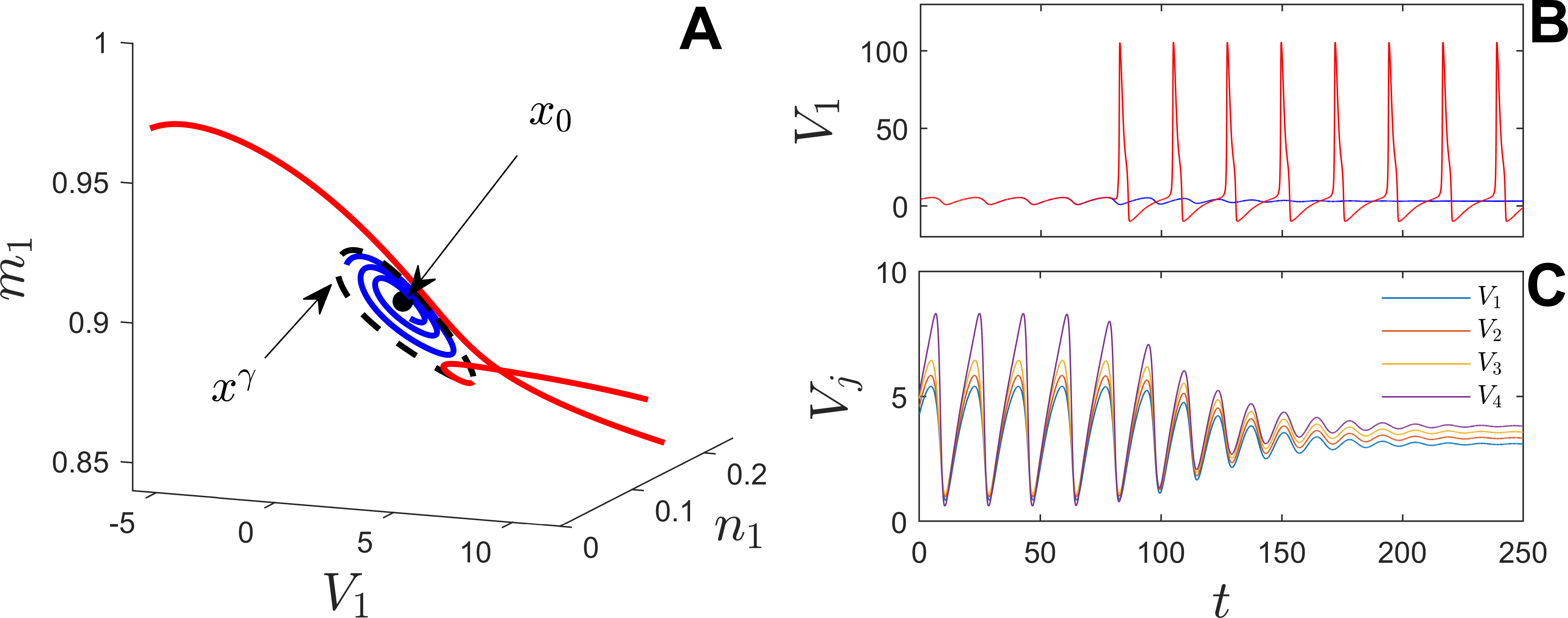}
\end{center}
\caption{ For the parameters used here, the model \eqref{hheqs} is close to a subcritical Hopf bifurcation.  An unstable periodic orbit, $x^\gamma$, separates solutions that decay to the stable fixed point $x_0$ and those that produce periodic action potentials.  The blue line is on the intersection of the unstable manifold of the unstable periodic orbit and the stable manifold of the stable fixed point.  Panel B shows $V_1$ for the blue and red trajectores from panel A over time.  Panel C shows the the transmembrane voltage of each neuron associated with the blue trajectory from panel A.}
\label{hhplots}
\end{figure}

For the parameters used here, the model \eqref{hheqs} is near a subcritical Hopf bifurcation.  Panel A of Figure \ref{hhplots} shows the relevant topological features of this model, plotted for a subset of the state variables.  An unstable periodic orbit, $x^\gamma$, separates solutions that fire periodic action potentials from those that converge towards a stable fixed point, $x_0$.  In panel A, two initial conditions are chosen on the unstable manifold of the unstable periodic orbit (dashed black line).  The red line is in the basin of attraction of the stable periodic orbit.  The blue line is in the basin of attraction of the stable fixed point.  The slowest decaying eigenvalues of the fixed point are $\lambda_{1,2} = -0.035 \pm 0.461$, with $\lambda_{3,4} = -0.473 \pm 0.408$ being the next slowest.     Blue and red lines in panel B of Figure \ref{hhplots} show the time course of $V_1$ for the blue and red trajectories shown in panel A.  Panel C shows the the transmembrane voltage of each neuron associated with the blue trajectory from panel A.


\begin{figure}[htb]
\begin{center}
\includegraphics[height=2.0 in]{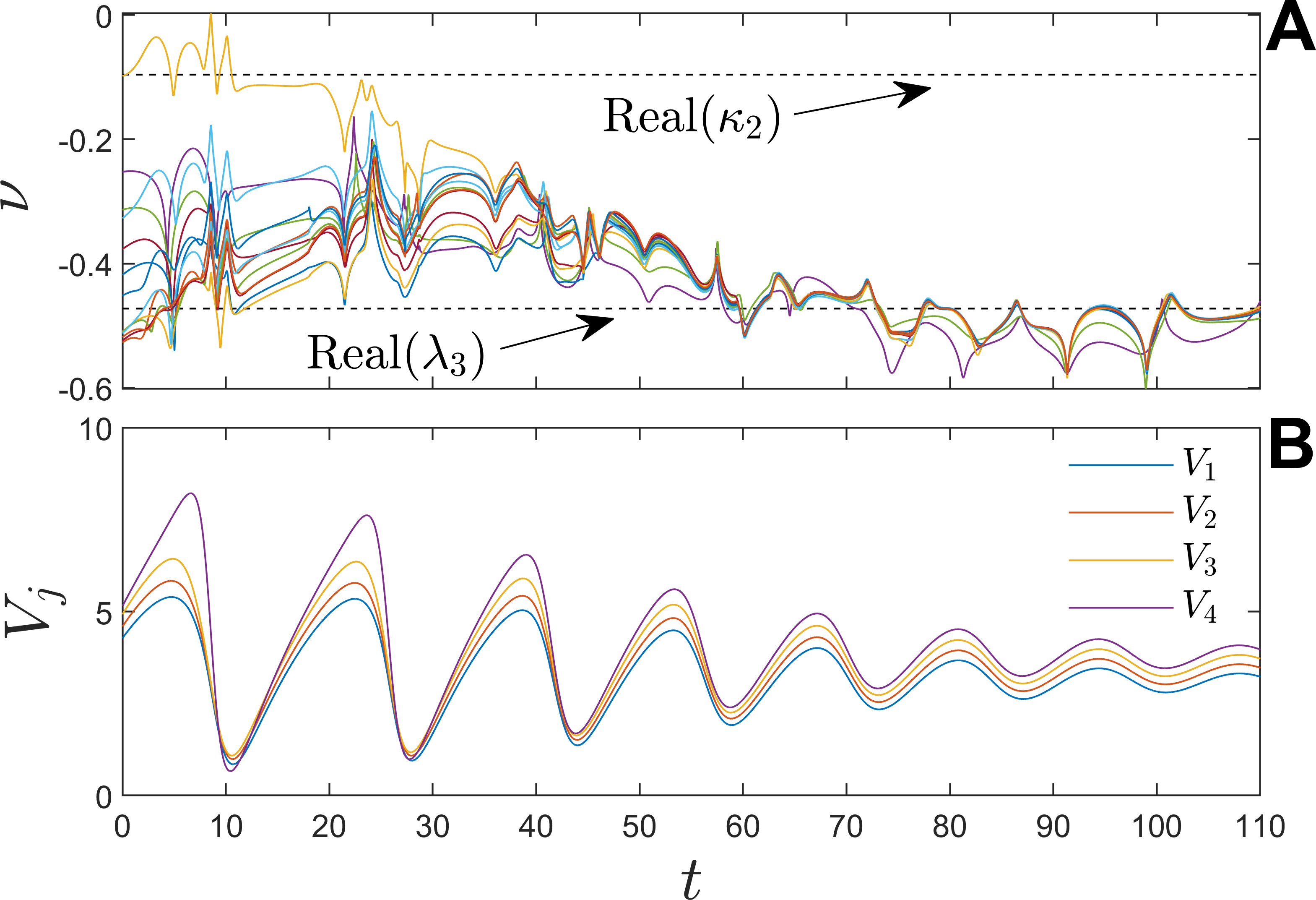}
\end{center}
\caption{For a representative sample trajectory on the slow manifold defined according to \eqref{slowstate2} panel A shows the decay rate $\nu$ associated with a set of perturbations transverse to the slow manifold.  Near the unstable periodic orbit, the slowest decaying modes decay at a rate close to $\kappa_2$ the slowest stable Floquet exponent.  Close to the fixed point, the modes decay at a rate close to ${\rm Real}(\lambda_3)$ which is associated with the third slowest decaying isostable coordinate.  Panel B shows the value of $V_1,\dots,V_4$ along this trajectory for reference.}
\label{convergencedata}
\end{figure}

The convergence rate of perturbations to the manifold defined according to \eqref{slowstate2} is considered along sample trajectories, with representative results shown in Figure \ref{convergencedata}.  For this slow manifold $\beta = 2$.  For a sample trajectory that starts near the unstable periodic orbit evolving on the slow manifold, $g_1$ and $g_2$ are computed according to \eqref{gper1} and \eqref{gper2}.  The orthogonal complement of the span of $g_1$ and $g_2$ is 14-dimensional.  Figure \ref{convergencedata} shows the decay rates for 14 perturbations that provide a basis for the orthogonal complement of $g_1$ and $g_2$.  For each, decay rates are computed according to \eqref{nuequation} during simulations of \eqref{xbardot}.  Colored traces in panel A of Figure \ref{convergencedata} show the moving average of each $\nu(t)$ computed over a 30 millisecond window.  For reference, colored lines in panel B show the transmembrane voltage of each neuron.  For states close to the periodic orbit, the slowest decaying modes have decay rates that are governed by the slowest stable Floquet exponent, $\kappa_2$.  As the system approaches the stable fixed point, the decay rates are governed by ${\rm Real}(\lambda_3)$ which is associated with the slowest decaying isostable coordinate that is truncated from the isostable coordinate reduction.  These decay rates are similar for all trajectories on the two-dimesional slow manifold.

\begin{figure}[htb]
\begin{center}
\includegraphics[height=2.7 in]{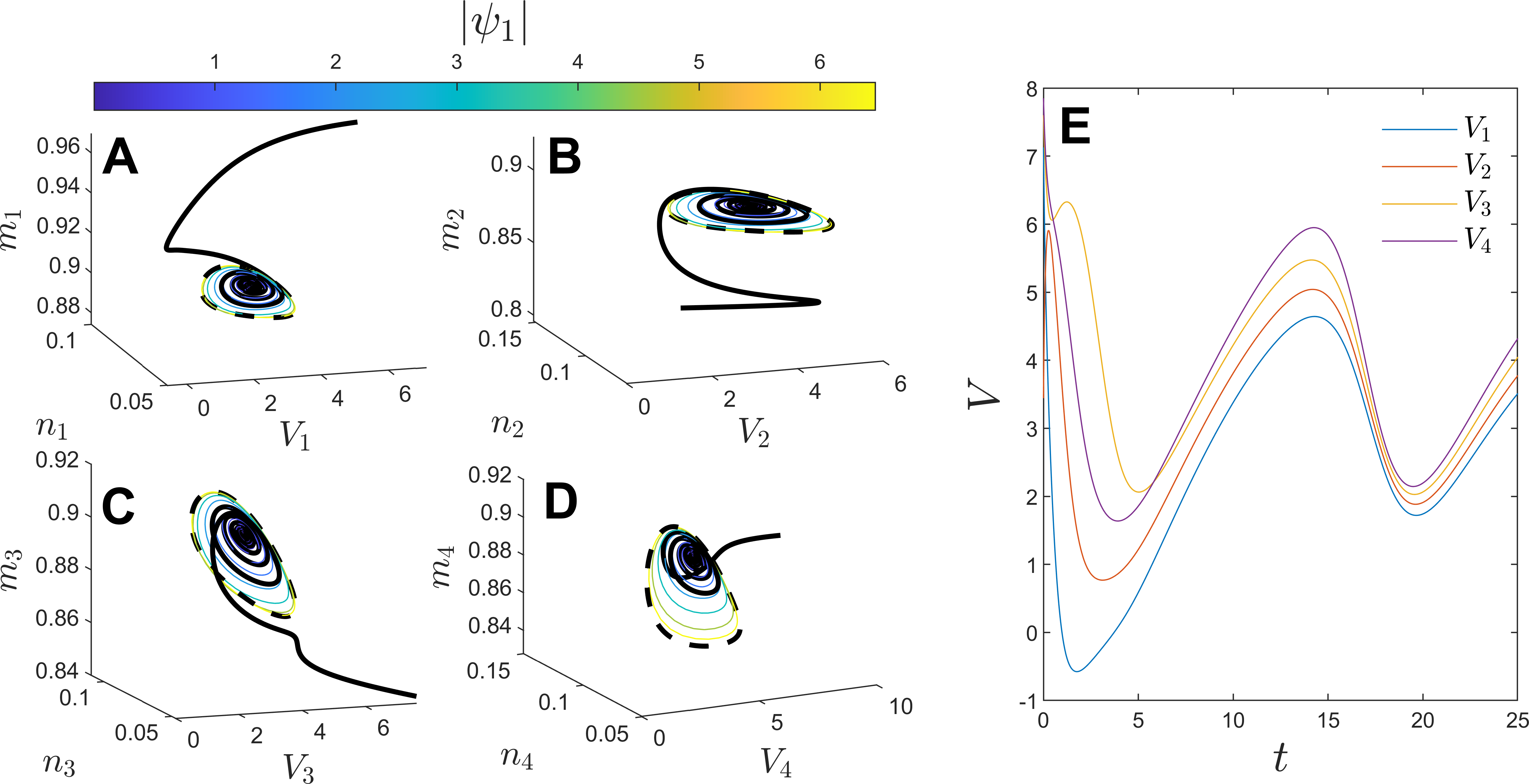}
\end{center}
\caption{Level sets of $|\psi_1|$  on the slow manifold of $x_0$ are shown as colored lines in panels A-D.  An initial condition in the basin of attraction of $x_0$ quickly converges to the slow manifold (black line).  For the same trajectory, panel E shows voltage traces for each neuron for this trajectory.}
\label{hhconvergence}
\end{figure}

A set of initial conditions close to $x^\gamma$ on its unstable manifold are integrated forward in time.  The resulting trajectories are used to trace out the slow manifold of the fixed point $x_0$.  $\psi_1$ is  subsequently computed and level sets of $|\psi_1|$ are plotted as colored lines in panels A-D of Figure \ref{hhconvergence}.  For an initial condition in the basin of attraction of $x_0$, the black line in Panels A-D show the evolution of $V_j$, $n_j$, and $m_J$ for $j = 1,2,3,4$  illustrating rapid convergence to the slow manifold.  For this same initial condition, panel E shows the transmembrane voltage of each neuron over time.  Despite having substantially different initial conditions, the plot quickly begins to resemble the plot from panel C of Figure \ref{hhplots}, i.e.,~shown for a trajectory on the slow manifold.

The slow manifold of the fixed point $x_0$ can be used for model reduction by considering the evolution of the 16-dimensional model \eqref{hheqs} on the 2-dimensional slow manifold.  As an illustration of this point, we consider a control strategy designed to prevent \eqref{hheqs} from firing an action potential in response to an external sine wave input.  Let the input $u$ be comprised of a sinusoidal forcing with an additional control input
\begin{equation} \label{addinput}
    u(t) = \zeta \sin(\omega_f t) + i(t),
\end{equation}
where $\zeta$ and $\omega_f$ is the magnitude and frequency of the sinusoidal forcing, and $i(t)$ is an additional input that will be designed to control the state back to the basin of attraction of $x_0$ if the sinusoidal forcing causes it to leave.  To implement such a control strategy, we focus on the dynamics of the isostable coordinates on the slow manifold.  For an initial condition on the slow manifold, neglecting the faster decaying isostable coordinates, the slow isostable coordinate $\psi_1$ evolves according to a simplified version of \eqref{isodyn}
\begin{align} \label{isodynexamp}
    \dot{\psi_1} &=  \lambda_1 \psi_1 + I_{\psi}(\psi_1) u(t),
\end{align}
where $I_\psi(\psi_1) = \frac{\partial \psi_1}{\partial V_1} + \frac{\partial \psi_1}{\partial V_2}$ with partial derivatives evaluated at $x(\psi_1)$ (recall that input is only applied in the voltage equations of the first two neurons).  Note that while there are two slow isostable coordinates, $\psi_1 = \psi_2^*$ so that $I_\psi$ can be written as a function of only $\psi_1$ when restricted to the slow manifold.  $I_1(\psi_1)$ can be obtained by first computing a set of trajectories on the slow manifold (for instance, the blue trajectory from panel A of Figure \ref{hhplots} could be one such trajectory), computing $I_1$ along each trajectory using \eqref{isoeq} and using the resulting information to interpolate $I_\psi(\psi)$ during simulations of \eqref{isodynexamp}.  Equation \eqref{isodyn} cannot be used to consider states beyond the basin of attraction of $x_0$, however, since state approaches the unstable periodic orbit $x^\gamma$ as $|\psi_1|$ approaches infinity.  In order to consider states beyond the basin of attraction, once $|\psi_1|$ increases beyond some threshold, a phase reduction of the form \eqref{isophasedyn1} and \eqref{isophasedyn2} will be considered
\begin{align} \label{phasedynexamp}
    \dot{\theta} &= \omega  + z(\theta) u(t), \nonumber \\
    \dot{\hat{ \psi }} &= \kappa \hat{\psi} + I_{\hat{\psi}} (\theta) u(t),  
\end{align}
where $\theta$ is the phase associated with the unstable periodic orbit, $\hat{\psi}$ is the isostable coordinate associated with its unstable Floquet multiplier $\kappa = 0.15$, and $\omega = 0.347$ rad/ms is the natural frequency.  Once again, because input is only applied to the voltage equations of the first two neurons $z(\theta) = \frac{\partial \theta}{\partial V_1} + \frac{\partial \theta}{\partial V_2}$ and  $I_{\hat{\psi}}(\theta) = \frac{\partial \hat{\psi}}{\partial V_1} + \frac{\partial \hat{\psi}}{\partial V_2}$ with partial derivatives evaluated at $x^\gamma(\theta)$.  Note that  Equation \eqref{phasedynexamp} is only valid for states close to the unstable periodic orbit, meaning that $\hat{\psi}$ must be small.  In \eqref{phasedynexamp}, $z(\theta)$ and $I_{\hat{\psi}}$ can be obtained by finding $Z(\theta)$ and $I_1(\theta)$ along $x^\gamma(t)$ using   \eqref{zidotper} and extracting the appropriate components.

The overall reduced order model for \eqref{hheqs} is implemented as follows:~Equation \eqref{isodynexamp} is used when $|\psi_1|\leq 6$, i.e.,~ when the state is far from the unstable periodic orbit.  If $|\psi_1|$ becomes larger than 6, the state is close to $x^\gamma$ and the state variables of \eqref{isodynexamp} are converted to those of the phase reduction \eqref{phasedynexamp}.  For this reduction, $\hat{\psi} < 0$, corresponds to states in the basin of attraction of $x_0$ and $\hat{\psi} \geq 0$ are outside the basin of attraction.  When $\hat{\psi}$ drops beyond a predetermined threshold, the state variables of \eqref{phasedynexamp} are converted back to those of isostable reduction \eqref{isodynexamp}.  Having determined an appropriate reduced order model, a relatively simple control input for preventing action potentials is as follows:~for an initial condition  in the basin of attraction of $x_0$, set the controller to be off so that $i(t) = 0$;  when Equation \eqref{phasedynexamp} is active and $\hat{\psi} > \psi_{\rm on}$, turn the controller on; when the controller is on let
\begin{equation} \label{controlhh}
   i(t) = \begin{cases} -C {\rm sign}(I_1(\psi_1) \psi_1^* + I_1^*(\psi_1) \psi_1), & \text{when Equation} \; \eqref{isodynexamp} \text{ is active}, \\
    -C {\rm sign}( I_0 ), & \text{when Equation} \; \eqref{phasedynexamp} \text{ is active},
    \end{cases}
\end{equation}
where $C$ is the applied control effort and $^*$ denotes the complex conjugate;
after the system returns to the basin of attraction and Equation \eqref{isodynexamp} is active, if $|\psi_1| < \psi_{\rm off}$, turn the controller off.  The intuition behind the controller \eqref{controlhh} is straightforward:~when \eqref{phasedynexamp} is active, the controller drives $\hat{\psi}$ to smaller values to bring it back within the basin of attraction.  When \eqref{isodynexamp} is active, noticing that $\frac{d\psi_1}{dt}|\psi_1|^2 = \frac{d}{dt} (\psi_1 \psi_1^*) = (\lambda_1 + \lambda_1^*)| \psi_1| + u (I_1(\psi_1) \psi_1^* + I_1^*(\psi_1)\psi_1)$, the controller drives $|\psi_1|$ to smaller values, bringing it further away from the boundary of the basin of attraction.

\begin{figure}[htb]
\begin{center}
\includegraphics[height=2.1 in]{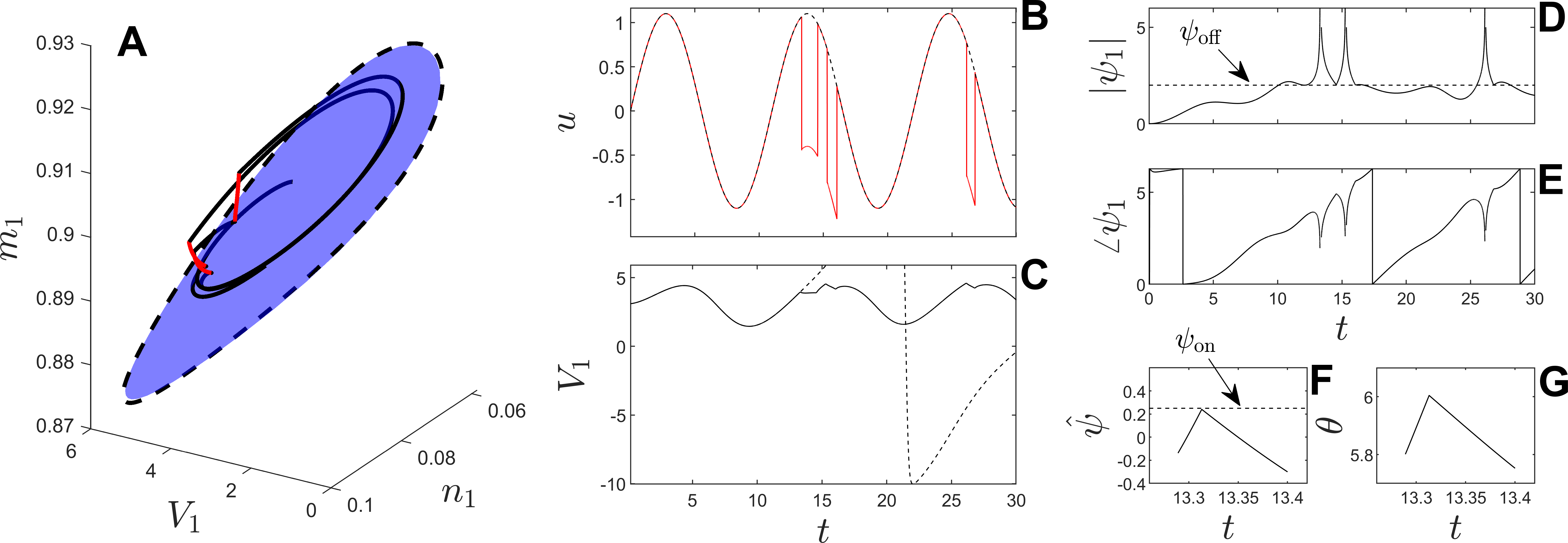}
\end{center}
\caption{Illustration of a control strategy for preventing action potentials in the model \eqref{hheqs} as described in the text.  Panel A shows the slow manifold for a subset of the state variables as a blue surface. The unstable periodic orbit (dashed line) is the boundary of the basin of attraction on the slow manifold.  The solid line shows the state in response to the control input.  In panel B, a sine wave is shown in black.  The red line shows $u(t)$ which differs from the sine wave when the control is active.  Panel C provides a plot of $V_1$ during the course of the simulation.  The dashed lines shows a comparison simulation taking $i(t)=0$ (only the sine wave is applied) with a resulting action potential occurring at approximately $t = 20$.  Panels D-G show the reduced order coordinates throughout the simulation with a detailed interpretation given in the text.}
\label{controlfig}
\end{figure}

Figure \ref{controlfig} illustrates the aforementioned control strategy. In panel A, the shaded blue region represents the states that are handled with the reduction \eqref{isodynexamp}.  It comprises most of the slow manifold, except for a small sliver near unstable periodic orbit (black dashed line).  All other states are handled by the phase-amplitude reduction \eqref{phasedynexamp}.  For the sinusoidal input from \eqref{addinput}, $\zeta = 1.1$ and $\omega_f = 0.57$.  Over a 30 ms simulation, panel B shows the sine wave in black, along with the value of $u$ as defined in \eqref{addinput} in red.  The sudden drops in the applied control represent moments during which the controller is on, bringing the state back within the basin of attraction.  The solid line in panel A show the trajectory in response to the applied input.  It is colored red when the control is on and black when the control is off.  Panel C shows a plot of $V_1$ over time for this simulation.  The dashed line shows the same simulation when the control is off for the entire time (i.e.,~just the sine wave is applied) illustrating that the sinusoidal input causes an action potential.  Panels D-G show the reduced order coordinates throughout the simulation.  Panels D and E show $|\psi_1|$ and the argument of $\psi_1$ when Equation \eqref{isodynexamp} is active.  Panels F and G show representative plots of the reduced order coordinates when \eqref{phasedynexamp} is active.  Note that \eqref{isodynexamp} and \eqref{phasedynexamp} are never simultaneously active which explains the gaps in each graph.

As an interpretation of the data in Panels D-G of Figure \ref{controlfig}, for the first 13 milliseconds the state stays  within the basin of attraction of the fixed point, and no control is applied.  At approximately 13.29 seconds, the state approaches the unstable periodic orbit with $|\psi_1|$ increasing beyond the threshold to switch from \eqref{isodynexamp} to \eqref{phasedynexamp}.    The sine wave input causes the state to exit the basin of attraction at approximately $\hat{\psi} = 13.3$ ms when $\hat{\psi} = 0$.  Once  $\hat{\psi}$ increases beyond $\psi_{\rm on} = 0.25$, the controller turns on, driving the state back inside the basin of attraction.   $\hat{\psi}$ decreasing below the threshold $-0.3$ prompts a switch from \eqref{phasedynexamp} to \eqref{isodynexamp}.  The control remains on until approximately $t = 14.6$ ms when $|\psi_1|$ crosses $\psi_{\rm off} = 2$.  This process repeats throughout the remainder of the simulation, with the control turning on as necessary to prevent the system from firing an action potential.  Note that for this simulation, the state is not prevented from leaving the basin of attraction of the fixed point, however, this could be mandated by setting $\psi_{\rm on} = 0$ and rerunning the simulation.

\subsection{Controlling an Oscillating Circadian Model to its Phaseless Set} \label{circsec}

\begin{figure}[htb]
\begin{center}
\includegraphics[height=2.7 in]{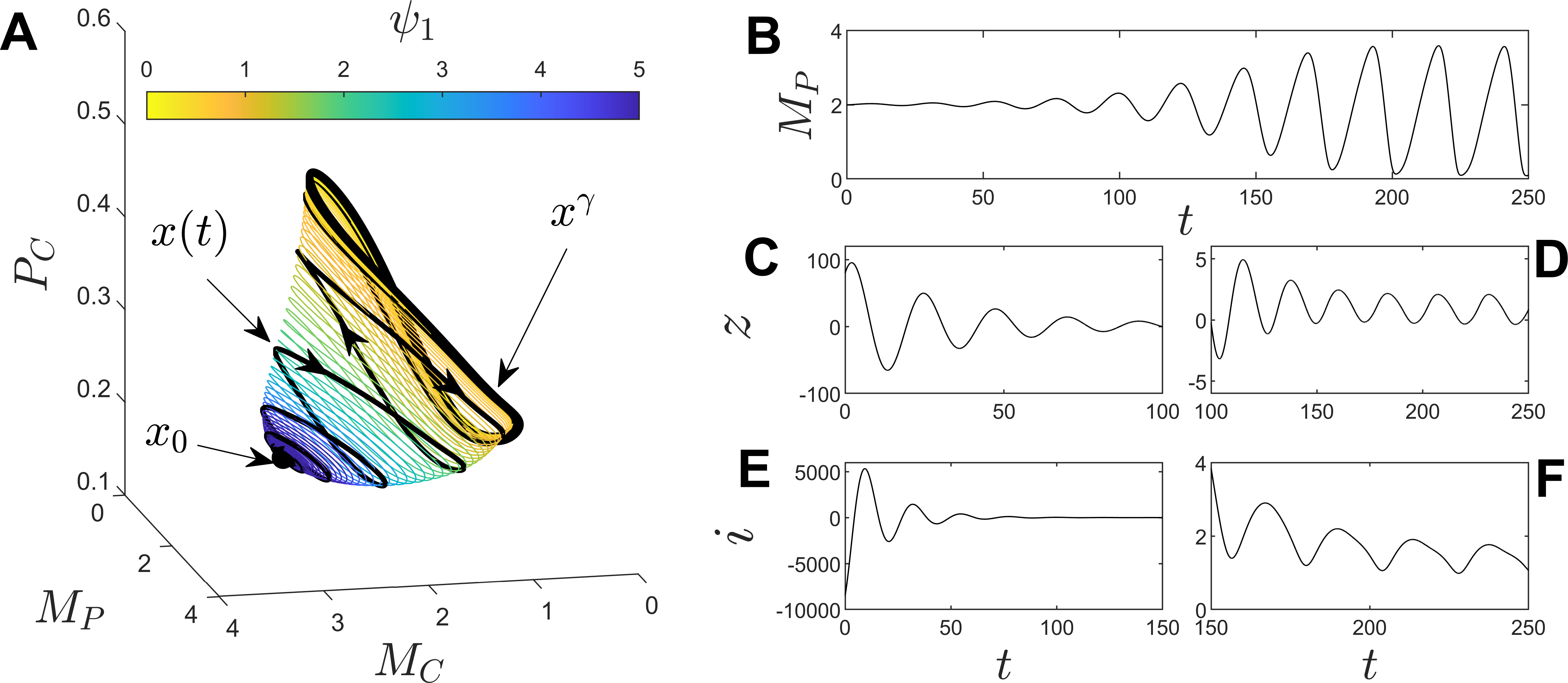}
\end{center}
\caption{In panel A, the trajectory $x(t)$ is on this slow manifold, starting near $x_0$ and ending near $x^\gamma$.  The colored lines show level sets of $\psi_1$ on the slow manifold.  Panel B shows the variable $M_p$ plotted along the trajectory $x(t)$. Recalling that an additive control input is added to $M_P$, Panels C-F show $z \equiv \frac{\partial \theta}{\partial M_P}$ and $i \equiv \frac{\partial \psi_1}{\partial M_P}$ plotted along the trajectory $x(t)$. These terms are relevant for reduced order modeling of the dynamics on the slow manifold. }
\label{circslowmanifold}
\end{figure}

Next, consider a 16 variable model \eqref{a16} that characterize the dynamical behavior of regulatory loops that govern the Per, Cry, Bmal1, and Clock genes \cite{lelo03} that give rise to circadian rhythms.  Full model equations are provided in Appendix \ref{circapx}.  Time is in units of hours.  An additive control input $u(t) \geq 0$ is applied to the $M_B$ concentration dynamics.  For the parameters used here, when $u(t) = 0$ this model has an unstable fixed point, $x_0$, shown in Panel A of Figure \ref{circslowmanifold} with unstable eigenvalues $\lambda_{1,2} =  0.0295 \pm 0.2776i$.  An initial condition on the unstable manifold of $x_0$ converges to a stable periodic orbit $x^\gamma$ in forward time.    This periodic orbit has a natural frequency of 0.262 radians/hour with a period of $T = 23.96$ hours.      The slowest decaying isostable coordinates are $\kappa_1 = -0.029$, $\kappa_2 = -0.065$, and $\kappa_3 = -0.1238$.  As compared to the example from Section \ref{hhsec} the gap between $\kappa_1$ and the faster Floquet exponents is not as pronounced.  Nevertheless the slowest Floquet exponent $\kappa_1$ can still be used to define a slow decaying isostable coordinate $\psi_1$ according to \eqref{iso1def} and an associated slow manifold on which $\psi_k = 0$ for $k \geq 2$.  

\begin{figure}[htb]
\begin{center}
\includegraphics[height=2.7 in]{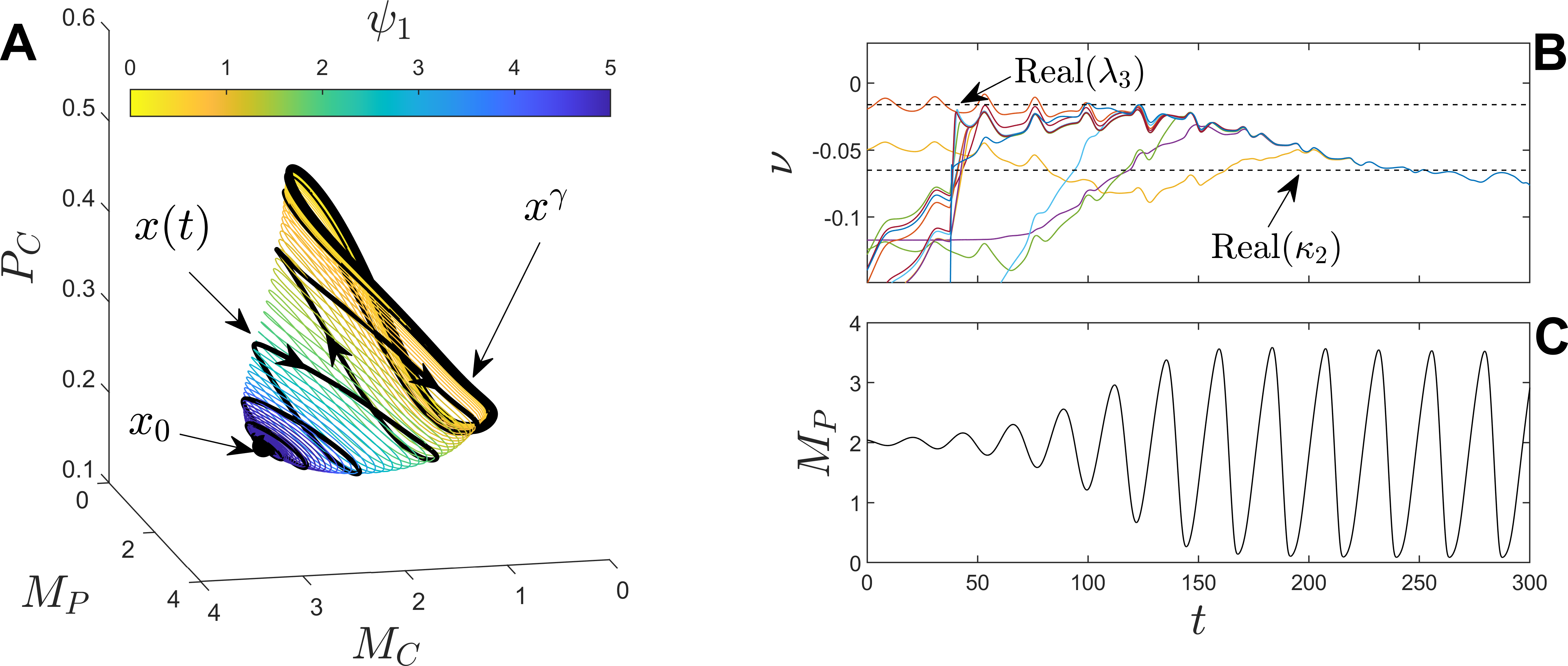}
\end{center}
\caption{A slow manifold of \eqref{a16} is given by the intersection of the unstable manifold of the unstable fixed point, $x_0$, and the stable manifold of the stable periodic orbit, $x^\gamma$. In panel A, the trajectory $x(t)$ is on this slow manifold, starting near $x_0$ and ending near $x^\gamma$.  The colored lines show level sets of $\psi_1$ (defined in reference to the stable periodic orbit) on the slow manifold.  Panel B shows the decay rate $\nu$ associated with a set of perturbations transverse to the slow manifold.  Panel C shows the value of $M_P$ along this trajectory for reference.}
\label{circdecayrates}
\end{figure}

Panel A of Figure \ref{circdecayrates} shows the evolution of a trajectory, $x(t)$ on the unstable manifold of $x_0$.  Level sets of $\psi_1$, (computed in reference to the stable periodic orbit) are shown for reference on the unstable manifold.   The dimension of the unstable manifold is $\beta = 2$.  For a sample trajectory that starts near the unstable fixed point evolving on the slow manifold, $g_1$ and $g_2$ are computed according to \eqref{findg}.  The orthogonal complement to the span of $g_1$ and $g_2$ is 14-dimensional.  The convergence rate of perturbations to the unstable manifold  is computed according to \eqref{nuequation} during simulations of \eqref{xbardot} along sample trajectories, with representative results shown in Figure \ref{convergencedata}.  Panel B of Figure \ref{circdecayrates} shows the moving average over a 100 hour window of decay rates for 14 perturbations that provide a basis for the orthogonal complement of $g_1$ and $g_2$.  These decay rates are computed according to \eqref{nuequation} during simulations of \eqref{nuequation}.  Panel C shows the value of $M_P$ on this trajectory for reference.  For states close to the unstable fixed point, the slowest decaying modes have decay rates that are governed by ${\rm Real}(\lambda_3)$, which is the slowest stable eigenmode.  As the system approaches the stable periodic orbit, the decay rates are governed by ${\rm Real}(\kappa_2)$ which is associated with the slowest decaying isostable coordinate that is truncated from the isostable coordinate reduction.  These decay rates are similar for all trajectories on the two-dimensional slow manifold.

 The slow manifold described above is used to obtain a reduced order model to capture the response to the perturbation $u(t)$.  An example trajectory on the slow manifold is shown in panel A of Figure \ref{circslowmanifold} with $M_p$ given along this trajectory for reference.   Panels C-F show relevant information associated with the trajectory $x(t)$.   Recall that the an additive control input is applied to the $M_B$ dynamics; as such, for model reduction purposes $z \equiv  \frac{\partial \theta}{\partial M_B}$ and $i \equiv \frac{\partial \psi_1}{\partial M_B}$ are the relevant curves.  Panels C and D (resp.~E and F) show $z$ (resp.,~i) along the trajectory $x(t)$.  Note that $i$ is strictly positive in the latter portions of this trajectory (near the periodic orbit) but takes both positive and negative values in the earlier portions (near the fixed point).  This feature will become important when discussing the control objective of driving the system from the periodic orbit to the unstable fixed point.

The slow manifold of the periodic orbit $x^\gamma$ can be used for reduced order modeling by considering the evolution of the 16 dimensional model \eqref{a16} on the 2-dimensional slow manifold.  The dynamics of the phase and isostable coordinates on the slow manifold follow the general form \eqref{isophasedyn1} and \eqref{isophasedyn2}
\begin{align} \label{phasereduced}
    \dot{\theta} &= \omega + z(\theta,\psi_1) u(t), \nonumber \\
    \dot{\psi_1} &= \kappa_1 \psi_1 + i(\theta,\psi_1) u(t).
\end{align}
Above, $z(\theta,\psi_1) = \frac{\partial \theta}{\partial M_P}$ and $i(\theta,\psi_1) = \frac{\partial \psi_1}{\partial M_P}$, where the partial derivatives are evaluated at $x(\theta,\psi_1)$ on the slow manifold.   Both $z(\theta,\psi_1)$ and $i(\theta,\psi_1)$ can be computed by obtaining a set of trajectories on the slow manifold (for instance, the blue trajectory from panel A of Figure \ref{circslowmanifold} could be one such trajectory), computing $I_1$ and $Z$ along each trajectory using \eqref{zidotper} and using the resulting information to interpolate $z$ and $i$ to perform simulations of \eqref{phasereduced}.

\begin{figure}[htb]
\begin{center}
\includegraphics[height=2.7 in]{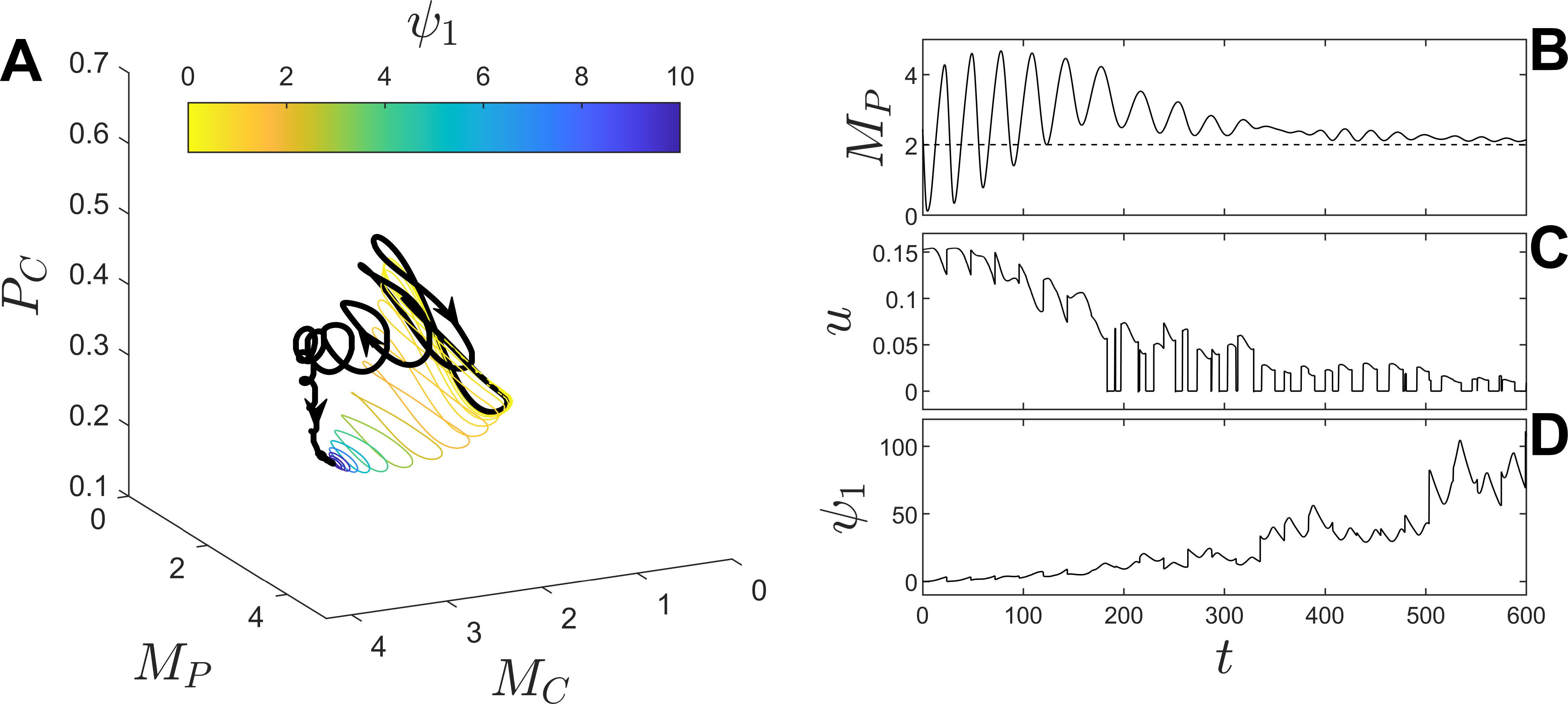}
\end{center}
\caption{  Applying the control strategy \eqref{circcontrol} to the circadian model \eqref{a16}.  In panel A, the controlled trajectory is plotted along side level sets of $\psi_1$ on the slow manifold.  The slow manifold is not as strongly attracting as in the previous example, but nonetheless provides a reduced order representation that is accurate enough to achieve the control objective of driving the system to the unstable fixed point.  Panels B, C, and D show $M_P$, $u$, and $\psi_1$ over time during this simulation.  The value of $\psi_1$ from \eqref{phasereduced} is estimated and updated every $T$ time units from measurements of $P_C(t)$ over the previous cycle.    }
\label{circcontrol_figure}
\end{figure}

The control objective for this example is to drive the state of \eqref{a16} from an initial condition on the periodic orbit to $x_0$.  For oscillatory systems, the fixed point is sometimes called the phaseless set, as all isochrons converge to this point \cite{winf01}, \cite{osin10}.  On the slow manifold, all trajectories converge to $x_0$ in the limit as $\psi_1$ approaches infinity.  As such, using the reduced order model \eqref{phasereduced}, driving the system from $x^\gamma$ to $x_0$ is equivalent to increasing $\psi_1$ as much as possible after starting from $\psi_1 = 0$.  From the level sets of $\psi_1$ in panel A of Figure \ref{circslowmanifold}, when $\psi_1>5$, the $x(\theta,\psi_1)$ is very close to $x_0$.  An additional constraint will be added requiring $u(t) \geq 0$ since it is generally easier to increase a given mRNA concentration (by adding more), than to decrease it.  With the aforementioned points in mind, the specific control input used is
\begin{equation} \label{circcontrol}
    u(t) =  \begin{cases}\frac{0.4 \omega}{M(\psi_1)} , & \text{if } {\rm sign}(i(\theta,\psi_1)) > 0, \\
    0, & \text{otherwise},
    \end{cases}
\end{equation}
where $M(\psi_1) = \max_\theta (|z(\theta,\psi_1)|)$.  Intuitively, when  ${\rm sign}(i(\theta,\psi_1)) > 0$, a positive input is applied to drive the system to larger isostable coordinates (closer to $x_0$).  The magnitude of the input is limited so that $z(\theta,\psi_1) u(t)$ cannot exceed 40 percent of $\omega$ and will shrink as the state gets closer to the target.   Note that when applying this control strategy to the full equations  \eqref{a16}, $\theta$ and $\psi_1$ are not directly accessible.  As such, the reduced order model \eqref{phasereduced} is simulated alongside \eqref{a16}, with the values of $\theta(t)$ and $\psi_1(t)$ used to determine the applied control.  Because \eqref{phasereduced} does not perfectly approximate the phase and isostable coordinates of the full order simulations, once every $T$ units of simulation time, the amplitude and oscillation timing of $P_C(t)$ over the previous cycle is used to update $\theta$ and $\psi_1$ in \eqref{phasereduced}.  Figure \ref{circcontrol_figure} provides a representative example of this control strategy applied to the full model \eqref{a16}.  Panel C shows the applied input and panel D shows $\psi_1$ over the course of the simulation.  Recall that $\psi_1$ is estimated from data every $T$ time units which explains the vertical lines in panel D.  Initially when the state is close to $x^\gamma$, $u$ is positive which is consistent with the fact that $i(\theta,\psi_1) > 0$ when $\psi_1$ is near zero.  Eventually as $\psi_1$ increases, $i(\theta,\psi_1) < 0$ in some places, which is consistent with the fact that $u = 0$ at some times as the simulation progresses.  Panel A shows the trajectory during the simulation, superimposed on levels sets of $\psi_1$ for states on the slow manifold.  Compared with the previous example from Section \ref{hhsec} the state does not converge as fast to the slow manifold due to smaller gap between $\kappa_1$ and $\kappa_2$.  Nonetheless, the reduced order model still enables a control strategy that can successfully achieve the control objective.  Panel B shows $M_p$ over time with the dashed horizontal line indicating the value of $M_P$ at the unstable fixed point for reference.

\begin{figure}[htb]
\begin{center}
\includegraphics[height=2.5 in]{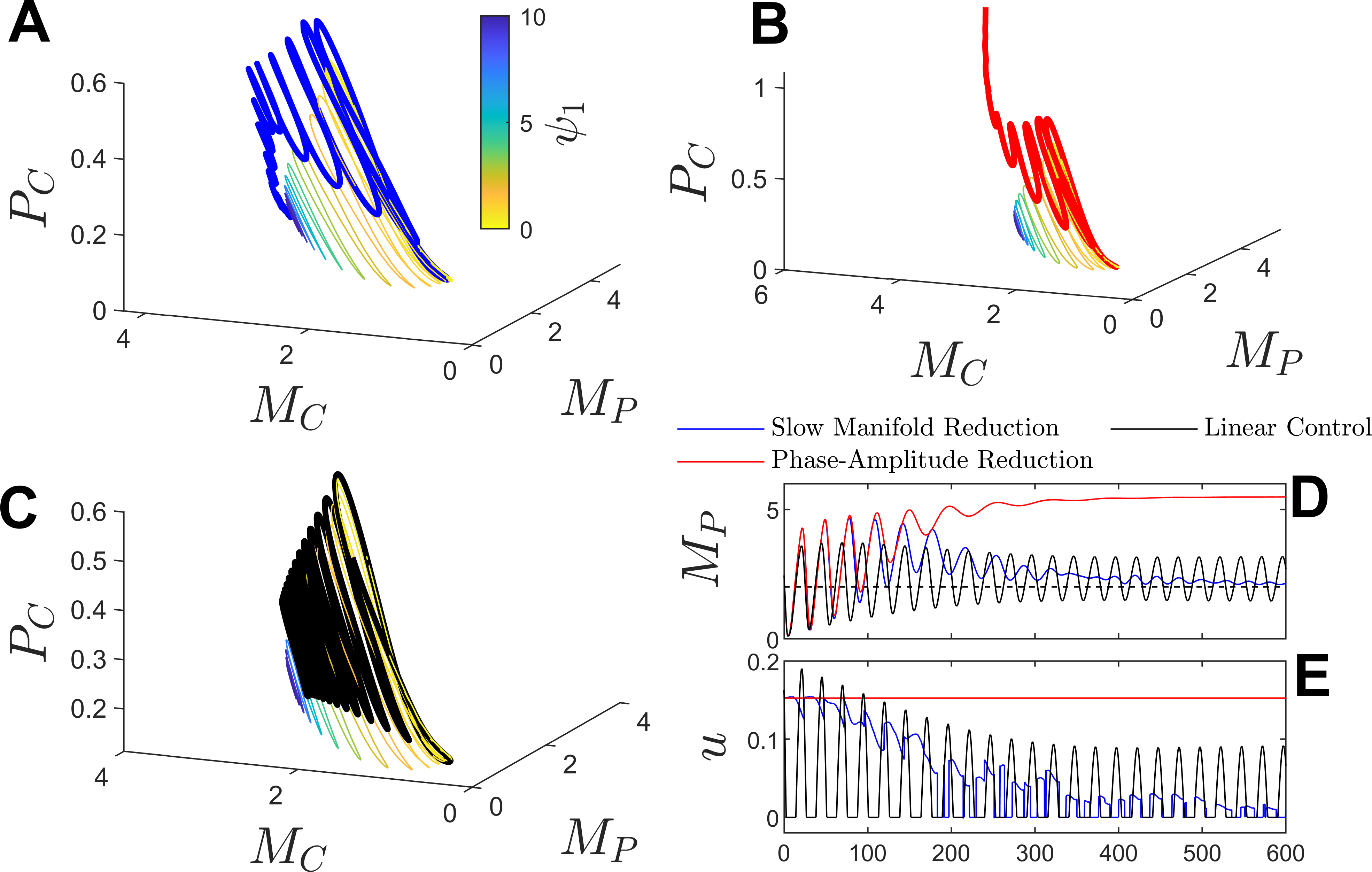}
\end{center}
\caption{ Comparison of two alternative control strategies for driving the state from the stable periodic orbit to the unstable fixed point (i.e.,~phaseless set).  Panels A,D, and E show the same results from figure \ref{circcontrol_figure} for comparison.  Panel B shows results when implementing a standard phase-amplitude reduction of the form \eqref{standardphase}.  Panel C shows the results of the control \eqref{lincontrol} based on a linearization of the model equations.  For each control strategy, panel D shows a plot of one $M_P$ over time, with the horizontal dashed line being the value at the unstable fixed point.  Panel E shows the control applied for each strategy.  Both the linear control strategy and the standard phase-amplitude-based control strategy do not achieve the control objective, necessitating the use of \eqref{phasereduced} as the reduced order model for this example.}
\label{differentcontrols}
\end{figure}

For comparison, two additional control strategies are considered in Figure \ref{differentcontrols}.  For the first, instead of evaluating the gradients of the phase and amplitude coordinates on the slow manifold, they are evaluated only on the periodic orbit yielding a model of the form
\begin{align} \label{standardphase}
    \dot{\theta} &= \omega + z(\theta,0) u(t), \nonumber \\
    \dot{\psi_1} &= \kappa_1 \psi_1 + i(\theta,0) u(t).
\end{align}
Equation \eqref{standardphase} is similar to reduced order models  used in \cite{wils16isos}, \cite{mong19b}, \cite{shoh21}.  While Equation \eqref{standardphase} is generally easier to implement than \eqref{phasereduced} because it requires less information about the system, it is only valid in a close neighborhood of the periodic orbit.  Using the reduced order model \eqref{standardphase}, a similar control strategy is applied to control the state of the full model to the unstable fixed point
\begin{equation} \label{circcontrolphaseonly}
    u(t) =  \begin{cases}\frac{0.4 \omega}{M(0)} , & \text{if } {\rm sign}(i(\theta,0)) > 0, \\
    0, & \text{otherwise},
    \end{cases}
\end{equation}
with results shown in panels B D and E Figure \ref{differentcontrols}.  Panel E shows the applied control (red line) is constant in time, which is consistent with the fact that $i(\theta,\phi)$ is positive for all $\theta$ when $\phi = 0$.  In panel B, the resulting trajectory is plotted along side level sets of $\psi_1$ on the slow manifold.  As compared to the results from Figure \ref{circcontrol_figure} (shown in panel A for reference) the controlled trajectory is initially similar, but ultimately misses the target indicating that knowledge of the behavior on the entire slow manifold is necessary to successfully implement this control strategy.  

A second control strategy based on a local linearization about $x_0$ is also considered.  The linearized model takes the form
\begin{equation}
    \Delta \dot{x} = A \Delta x + B u(t),
\end{equation}
where $\Delta x = x-x_0$, $A$ is the Jacobian evaluated at $x_0$, and $B \in \mathbb{R}^{16} = [ 0 \; 0\; 1 \; 0 \; \dots \; 0]^T$ (recall that the input is added to the third variable).   The unstable eigenvalues of the fixed point are $\lambda_{1,2} =  0.0295 \pm 0.2776i$.  Let $v_j$ and $w_j$ be the left and right eigenvectors associated with the eigenvalue $\lambda_j$ and let 
\begin{equation} \label{yest}
    \Upsilon  = [w_1 \; w_2]^T \Delta x \in \mathbb{C}^2,
\end{equation} 
be the coordinates associated with the unstable eigenmodes.  Changing variables, the dynamics of $\Upsilon$ follow
\begin{equation} 
    \dot{\Upsilon} = \begin{bmatrix} \lambda_1 & 0 \\ 0 & \lambda_2 \end{bmatrix} \Upsilon + \tilde{B} u(t),
\end{equation}
where $\tilde{B} = [w_1 \; w_2]^T B$.  A linear control strategy that is similar in spirit to \eqref{circcontrol} can then be designed with the goal of driving $\Upsilon$ to zero so that the system tends towards the fixed point.  Noting that the components of $\Upsilon = [\Upsilon_1 \; \Upsilon_2]^T$ are complex conjugates, a simple control scheme designed to drive $\Upsilon_1$ to lower magnitude values is
\begin{equation} \label{lincontrol}
u(t) =  \max( - \zeta \big[ \Upsilon_1^*(t) \tilde{B}_1 + \Upsilon_1(t) \tilde{B}^*_1 \big],0),
\end{equation}
where $\tilde{B}_1$ is the first entry of $\tilde{B}$, $^*$ denotes the complex conjugate, and $\zeta$ is a positive constant.  The control \eqref{lincontrol} seeks to decrease the magnitude of $|\psi_1|$, with an input magnitude that decreases as the state approaches the fixed point.  This linear control strategy \eqref{lincontrol}  ignores the effect of stable eigenmodes, similar to the control using the slow manifold from  \eqref{circcontrol}.  The linear control strategy \eqref{lincontrol} is applied to \eqref{a16}, estimating the value of $y$ according to \eqref{yest} assuming direct access to all state variables.  Results are shown in panels C-E of Figure \ref{differentcontrols} taking $\zeta = 0.08$.  As seen in panel C, the state starts to tend towards $x_0$ but does not reach it, owing to unmodeled nonlinearities in the system.  For the linear control strategy, larger values of $\zeta$ tend to destabilize the system, ultimately sending it far from the fixed point; smaller values of $\zeta$ yield steady state behavior with larger amplitude oscillations.

\section{Conclusion} \label{concsec}

This work investigates slow manifolds embedded in state space that are defined by the intersection of an unstable manifold of an unstable fixed point or periodic orbit and the stable manifold of a stable attractor.  Such a slow manifold can be approximated numerically by choosing a family of initial conditions near the unstable equilibrium, and integrating forward in time.  The dynamics of isostable coordinates (and phase coordinate as appropriate) on the slow manifold can be considered via the transformation \eqref{isodyn} when the stable attractor is a fixed point or via \eqref{isophasedyn1} and \eqref{isophasedyn2} when the stable attractor is a periodic orbit.  These transformations ultimately yield a reduced order model that describes the behavior of the system near the slow manifold.  Two biologically motivated systems are considered in this work to illustrate this reduced order modeling framework.

For the coupled population of Hodgkin-Huxley neurons considered in Section \ref{hhsec}, the slow manifold defined according to \eqref{slowstate1} gives an estimate of the basin of attraction of the system's fixed point.  Other Koopman-based strategies have been proposed for providing basin of attraction estimates for nonlinear dynamical systems \cite{will15}, \cite{maur16}, but these methods would be difficult to scale to high dimensional systems due to either the use gridding the state space or due to the need to compute high order partial derivatives -- the computational effort of such approaches increases exponentially with the dimension of the system.  By contrast, restricting attention to a low-dimensional slow manifold allows one to circumvent issues caused by high dimensionality of the underlying system.   The resulting reduced order model \eqref{isodynexamp} in conjunction with the phase-based model \eqref{phasedynexamp} enables a simple control strategy that is able to drive the system state back to within the basin of attraction after escaping due to the influence of an external input.   

For the circadian model considered in Section \ref{circsec}, the reduced order model giving the dynamics on the slow manifold allows for the design of a control algorithm to drive the state from a stable periodic orbit to an unstable fixed point (i.e.,~the phaseless set).  As shown in the results from Figure \ref{differentcontrols}, it is essential to understand the behavior along the entire slow manifold to achieve this control objective.   Standard phase-based reduction methods and strategies based on local linearization do not contain sufficient information about the system dynamics and control strategies based on these reductions do not achieve the control objective.

For a general system, in order to identify a slow manifold defined according to \eqref{slowstate1} or \eqref{slowstate2}, it is first necessary to identify stable and unstable fixed points and periodic orbits present in a given dynamical system.  This can pose a challenge, particularly in high dimensional dynamical systems, where the identification of unstable equilibira can be challenging.  In this work, this computation is accomplished using bifurcation continuation techniques, starting, for instance from a stable fixed point and following the solution as a parameter is changed.  

The proposed model order reduction techniques will not be applicable to all dynamical systems.  Foremost, a sufficient timescale separation between the decay rates of the fast and slow isostable coordinates is necessary to obtain an accurate reduced order system.  Larger timescale separation will produce more accurate reduced order models, however, moderate timescale difference can still yield accurate models (for instance, in Section \ref{circsec} the slow Floquet exponent is $\kappa_1 = -0.029$ with the next slowest Floquet exponent being $\kappa_2 = -0.065$).


This material is based upon the work supported by the National Science Foundation (NSF) under Grant No.~CMMI-2140527.

\begin{appendices}

\section{Hodgkin-Huxley Model Equations} \label{hhapx}
\renewcommand{\thetable}{A\arabic{table}}  
\renewcommand{\thefigure}{A\arabic{figure}} 
\renewcommand{\theequation}{A\arabic{equation}} 
\setcounter{equation}{0}
\setcounter{figure}{0}

The dynamical model used in Section \ref{hhsec} consists of 4 coupled Hodgkin-Huxley neurons \cite{hodg52}.  The equations are:
\begin{align} \label{hheqsapx}
    C_m  \dot{V}_j &= -g_{\rm Na} m_j^3 h_j (V_j - V_{\rm Na}) - g_{\rm K} n_j^4 (V_j - V_{\rm K}) - g_{\rm L}(V_j - V_{\rm L})  - g_C(V_j - \bar{V}) + I_j + \alpha_j u(t), \nonumber \\
    \dot{n}_j &= \alpha_n(1-n_j) - \beta_n n_j, \nonumber \\
    \dot{m}_j &= \alpha_m(1-m_j) - \beta_m m_j, \nonumber \\
    \dot{h}_j &= \alpha_h(1-h_j) - \beta_h h_j,
\end{align}
for $j = 1,\dots 4$, where $V_j$ is the transmembrane potential (in mV), and $n_j$, $m_j$, and $h_j \in [0,1]$ are gating variables of the $j^{\rm th}$ neuron.  $C_m = 1 \mu{\rm F}/{\rm cm}^2$ is the membrane capacitance, and $I_j = -4.318 + 0.75j$ $\mu {\rm A}/{\rm cm}^2$ is the baseline current of each neuron.  Coupling is electrotonic \cite{john95} with $g_C = 1.5$ and $\bar{V} = \frac{1}{4} \sum_{j = 1}^4 V_j$.  Here $u(t)$ is an input that represents a transmembrane current applied to only the first two neurons so that $\alpha_j = 1$ for $j = 1,2$ and equals zero otherwise.

Maximal membrane conductances are
\begin{align}
g_{\rm Na}= 120 \; {\rm mS}/{\rm cm}^2, \nonumber \\
g_{\rm K}=  36 \; {\rm mS}/{\rm cm}^2, \nonumber \\
g_{\rm L}=  1 \; {\rm mS}/{\rm cm}^2,  
\end{align}
and
\begin{align}
V_{\rm Na}= 115 \; {\rm mV}, \nonumber \\
V_{\rm K}=  -12 \; {\rm mV}, \nonumber \\
V_{\rm L}=  10.599 \; {\rm mV},
\end{align}
are the reversal potentials of the associated ion channels.  The rate constants are functions of the transmembrane voltage
\begin{align}
\alpha_n &= 0.01(10-V)/(\exp((10-V)/10)-1), \nonumber \\
\beta_n &= 0.125\exp(-V/80), \nonumber \\
\alpha_m &= 0.1(25-V)/(\exp((25-V)/10)-1), \nonumber \\
\beta_m &= 4\exp(-V/18), \nonumber \\
\alpha_h &= 0.7\exp(-V/20), \nonumber \\
\beta_h &= 1/(\exp((-V+30)/10)+1). 
\end{align}

\section{Circadian Model Equations} \label{circapx}
\renewcommand{\thetable}{B\arabic{table}}  
\renewcommand{\thefigure}{B\arabic{figure}} 
\renewcommand{\theequation}{B\arabic{equation}} 
\setcounter{equation}{0}
\setcounter{figure}{0}

The circadian oscillator model used in Section \ref{circsec} was published in \cite{lelo03}.  The version used here has 16 coupled ordinary differential equations.  The state variables are as follows:~concentrations of  Per, Cry, and Bmal1 mRNA are designated by $M_P$, $M_C$, and $M_B$, respectively;  phosphorylated (resp.,~nonphosphorylated)  Per and Cry proteins in cytosol are designated by $P_{CP}$ and $C_{CP}$ (resp.,~$P_C$ and $C_C$);  concentrations of Per-Cry complex in cytosol and nucleus are designated by $PC_C$, $PC_N$, $PC_{CP}$, and $PC_{NP}$;  concentrations of BMAL1 in cytosol and nucleus are designated by $B_C$, $B_{CP}$, $B_N$, and $B_{NP}$;  Inactive complex between Per-Cry and Clock-Bmal1 in the nucleus is designated by $I_N$.  Subscripts $C$, $N$, $CP$ and $NP$ denote cytosolic, nuclear, cytosolic phosphorylated, and nuclear phosphorylated forms, respectively.    The model equations are:
\begin{align}
\dot{M}_P &= v_{sP} \frac{B_N^n}{K_{AP}^n + B_N^n} - v_{mP} \frac{M_P}{K_{mP} + M_P} - k_{dmp}M_P,   \nonumber \\
\dot{M}_C &= v_{sC}\frac{B_N^n}{K^n_{AC} + B_N^n} - v_{mC} \frac{M_C}{K_{mC} + M_C} - k_{dmc} M_C, \nonumber \\
\dot{M}_B &= v_{sB}\frac{K_{IB}^m}{K_{IB}^m + B_N^m} - v_{mB} \frac{M_B}{K_{mB} + M_B} - k_{dmb} M_B + u(t),  \nonumber \\
\dot{P}_C &= k_{sP}M_P - V_{1P} \frac{P_C}{K_p + P_C} + V_{2P}\frac{P_{CP}}{K_{dp} + P_{CP}} + k_4 PC_C - k_3 P_CC_C - k_{dn}P_C, \nonumber \\
\dot{C}_C &= k_{sC}M_C - V_{1C} \frac{C_C}{K_p + C_C} + V_{2C}\frac{C_{CP}}{K_{dp} + C_{CP}} + k_4 PC_C - k_3 P_CC_C - k_{dnc}C_C, \nonumber \\
\dot{P}_{CP} &= V_{1P} \frac{P_C}{K_p + P_C} - V_{2P} \frac{P_{CP}}{K_{dp} + P_{CP}} - v_{dPC} \frac{P_{CP}}{K_d + P_{CP}} - k_{dn}P_{CP}, \nonumber \\
\dot{C}_{CP} &= V_{1C} \frac{C_C}{K_p + C_C} - V_{2C} \frac{C_{CP}}{K_{dp} + C_{CP}} - v_{dCC} \frac{C_{CP}}{K_d + C_{CP}} - k_{dn}C_{CP}, \nonumber \\
\dot{PC}_C &= - V_{1PC} \frac{PC_C}{K_p + PC_C} + V_{2PC}\frac{PC_{CP}}{K_{dp} + PC_{CP}} - k_4 PC_C + k_3 P_C C_C \nonumber  \nonumber \\
&\quad+ k_2 PC_N - k_1PC_C - k_{dn}PC_C, \nonumber \\
\dot{PC}_N &= - V_{3PC} \frac{PC_N}{K_p + PC_N} + V_{4PC}\frac{PC_{NP}}{K_{dp} + PC_{NP}} - k_2 PC_N + k_1 PC_C \nonumber \nonumber \\
&\quad- k_7B_N PC_N + k_8I_N - k_{dn}PC_N, \nonumber \\
\dot{PC}_{CP} &= V_{1PC} \frac{PC_C}{K_p + PC_C} - V_{2PC} \frac{PC_{CP}}{K_{dp} + PC_{CP}} - v_{dPCC} \frac{PC_{CP}}{K_d + PC_{CP}} - k_{dn}PC_{CP}, \nonumber \\
\dot{PC}_{NP} &=- V_{3PC} \frac{PC_N}{K_p + PC_N} - V_{4PC}\frac{PC_{NP}}{K_{dp} + PC_{NP}} - v_{dPCN} \frac{PC_{NP}}{K_d + PC_{NP}} - k_{dn}PC_{NP},  \nonumber \\
\dot{B}_C &= k_{sB} M_B - V_{1B} \frac{B_C}{K_p + B_C} + V_{2B} \frac{B_{CP}}{K_{dp} + B_{CP}} - k_5B_C + k_6B_N - k_{dn}B_C, \nonumber \\
\dot{B}_{CP} &=  V_{1B} \frac{B_C}{K_p + B_C} - V_{2B} \frac{B_{CP}}{K_{dp} + B_{CP}} - v_{dBC} \frac{B_{CP}}{K_d + B_{CP}} - k_{dn}B_{CP}, \nonumber \\
\dot{B}_{N} &=  -V_{3B} \frac{B_N}{K_p + B_N} + V_{4B} \frac{B_{NP}}{K_{dp} + B_{NP}} + k_5 B_C - k_6 B_N - k_7 B_N PC_N  \nonumber  \nonumber \\
& \quad + k_8 I_N - k_{dn}B_N, \nonumber  \\
\dot{B}_{NP} &=  V_{3B} \frac{B_N}{K_p + B_N} - V_{4B} \frac{B_{NP}}{K_{dp} + B_{NP}} - v_{dBN} \frac{B_{NP}}{K_d + B_{NP}} - k_{dn}B_{NP}, \nonumber \\
\dot{I}_N &= -k_8 I_N + k_7 B_N PC_N-v_{dIN} \frac{I_N}{K_d + I_N} - k_{dn}I_N  \label{a16}.
\end{align}
Basal values listed in Supplementary Table 1 of \cite{lelo03} are used with the exception of $k_1 = 0.58$ and $k_2 = 2.0$, which determine the dynamics of the nonphysphorylated Per-Cry complex.  Units of time are in hours.  The light intensity impacts the maximum rate of Per expression and is taken to be $v_{sP} = 1.35$.  A control input $u(t)$ is added to the $M_B$ concentration dynamics.

\end{appendices}


\begin{thebibliography}{10}

\bibitem{ahme23}
T.~Ahmed, A.~Sadovnik, and D.~Wilson.
\newblock Data-driven inference of low-order isostable-coordinate-based dynamical models using neural networks.
\newblock {\em Nonlinear Dynamics}, 111(3):2501--2519, 2023.

\bibitem{brun16}
S.~L. Brunton, B.~W. Brunton, J.~L. Proctor, and J.~N. Kutz.
\newblock Koopman invariant subspaces and finite linear representations of nonlinear dynamical systems for control.
\newblock {\em PloS One}, 11(2), 2016.

\bibitem{cene22}
M.~Cenedese, J.~Ax{\aa}s, B.~B{\"a}uerlein, K.~Avila, and G.~Haller.
\newblock Data-driven modeling and prediction of non-linearizable dynamics via spectral submanifolds.
\newblock {\em Nature Communications}, 13(1):872, 2022.

\bibitem{farj18}
S.~Farjami, V.~Kirk, and H.~M. Osinga.
\newblock Computing the stable manifold of a saddle slow manifold.
\newblock {\em SIAM Journal on Applied Dynamical Systems}, 17(1):350--379, 2018.

\bibitem{feni79}
N.~Fenichel.
\newblock Geometric singular perturbation theory for ordinary differential equations.
\newblock {\em Journal of differential equations}, 31(1):53--98, 1979.

\bibitem{foia882}
C.~Foias, M.~S. Jolly, I.~G. Kevrekidis, G.~R. Sell, and E.~S. Titi.
\newblock On the computation of inertial manifolds.
\newblock {\em Physics Letters A}, 131(7):433--436, 1988.

\bibitem{foia88}
C.~Foias, G.~R. Sell, and R.~Temam.
\newblock Inertial manifolds for nonlinear evolutionary equations.
\newblock {\em Journal of Differential Equations}, 73(2):309--353, 1988.

\bibitem{guck75}
J.~Guckenheimer.
\newblock Isochrons and phaseless sets.
\newblock {\em Journal of Mathematical Biology}, 1(3):259--273, 1975.

\bibitem{guck09}
J.~Guckenheimer and C.~Kuehn.
\newblock Computing slow manifolds of saddle type.
\newblock {\em SIAM Journal on Applied Dynamical Systems}, 8(3):854--879, 2009.

\bibitem{hall16}
G.~Haller and S.~Ponsioen.
\newblock Nonlinear normal modes and spectral submanifolds: existence, uniqueness and use in model reduction.
\newblock {\em Nonlinear Dynamics}, 86:1493--1534, 2016.

\bibitem{hodg52}
A.~L. Hodgkin and A.~F. Huxley.
\newblock A quantitative description of membrane current and its application to conduction and excitation in nerve.
\newblock {\em J. Physiol.}, 117:500--44, 1952.

\bibitem{john95}
D.~Johnston and S.~M.-S. Wu.
\newblock {\em Foundations of Cellular Neurophysiology}.
\newblock MIT Press, Cambridge, MA, 1995.

\bibitem{jord07}
D.~Jordan and P.~Smith.
\newblock {\em Nonlinear Ordinary Differential Equations:~An Introduction for Scientists and Engineers}, volume~10.
\newblock Oxford University Press, Oxford, 2007.

\bibitem{kais21}
E.~Kaiser, J.~N. Kutz, and S.~Brunton.
\newblock Data-driven discovery of {K}oopman eigenfunctions for control.
\newblock {\em Machine Learning: Science and Technology}, 2021.

\bibitem{kape99}
T.~J. Kaper.
\newblock Systems theory for singular perturbation problems.
\newblock In {\em Analyzing multiscale phenomena using singular perturbation methods}, pages 85--131. AMS, 1999.

\bibitem{kutz16}
J.~N. Kutz, S.~L. Brunton, B.~W. Brunton, and J.~L. Proctor.
\newblock {\em Dynamic mode decomposition: data-driven modeling of complex systems}.
\newblock Society for Industrial and Applied Mathematics, Philadelphia, PA, 2016.

\bibitem{kval21}
M.~D. Kvalheim and S.~Revzen.
\newblock Existence and uniqueness of global {K}oopman eigenfunctions for stable fixed points and periodic orbits.
\newblock {\em Physica D: Nonlinear Phenomena}, page 132959, 2021.

\bibitem{lan13}
Y.~Lan and I.~Mezi{\'c}.
\newblock Linearization in the large of nonlinear systems and {K}oopman operator spectrum.
\newblock {\em Physica D: Nonlinear Phenomena}, 242(1):42--53, 2013.

\bibitem{lelo03}
J.~C. Leloup and A.~Goldbeter.
\newblock Toward a detailed computational model for the mammalian circadian clock.
\newblock {\em Proceedings of the National Academy of Sciences}, 100(12):7051--7056, 2003.

\bibitem{lusc18}
B.~Lusch, J.~N. Kutz, and S.~L. Brunton.
\newblock Deep learning for universal linear embeddings of nonlinear dynamics.
\newblock {\em Nature Communications}, 9(1):1--10, 2018.

\bibitem{maur16}
A.~Mauroy and I.~Mezi{\'c}.
\newblock Global stability analysis using the eigenfunctions of the {K}oopman operator.
\newblock {\em IEEE Transactions on Automatic Control}, 61(11):3356--3369, 2016.

\bibitem{maur13}
A.~Mauroy, I.~Mezi{\'c}, and J.~Moehlis.
\newblock Isostables, isochrons, and {K}oopman spectrum for the action--angle representation of stable fixed point dynamics.
\newblock {\em Physica D: Nonlinear Phenomena}, 261:19--30, 2013.

\bibitem{mezi13}
I.~Mezi{\'c}.
\newblock Analysis of fluid flows via spectral properties of the {K}oopman operator.
\newblock {\em Annual Review of Fluid Mechanics}, 45:357--378, 2013.

\bibitem{mezi19}
I.~Mezi{\'c}.
\newblock Spectrum of the {K}oopman operator, spectral expansions in functional spaces, and state-space geometry.
\newblock {\em Journal of Nonlinear Science}, pages 1--55, 2019.

\bibitem{mezi20}
I.~Mezi{\'c}.
\newblock Spectrum of the {K}oopman operator, spectral expansions in functional spaces, and state-space geometry.
\newblock {\em Journal of Nonlinear Science}, 30(5):2091--2145, 2020.

\bibitem{mong19b}
B.~Monga and J.~Moehlis.
\newblock Optimal phase control of biological oscillators using augmented phase reduction.
\newblock {\em Biological Cybernetics}, 113(1-2):161--178, 2019.

\bibitem{osin10}
H.~M. Osinga and J.~Moehlis.
\newblock Continuation-based computation of global isochrons.
\newblock {\em SIAM Journal on Applied Dynamical Systems}, 9(4):1201--1228, 2010.

\bibitem{park24}
Y.~Park and D.~Wilson.
\newblock N-body oscillator interactions of higher-order coupling functions.
\newblock {\em SIAM Journal on Applied Dynamical Systems}, 23(2):1471--1503, 2024.

\bibitem{pons20}
S.~Ponsioen, S.~Jain, and G.~Haller.
\newblock Model reduction to spectral submanifolds and forced-response calculation in high-dimensional mechanical systems.
\newblock {\em Journal of Sound and Vibration}, 488:115640, 2020.

\bibitem{schm10}
P.~J. Schmid.
\newblock Dynamic mode decomposition of numerical and experimental data.
\newblock {\em Journal of Fluid Mechanics}, 656:5--28, 2010.

\bibitem{soot17}
A.~Sootla and A.~Mauroy.
\newblock Geometric properties of isostables and basins of attraction of monotone systems.
\newblock {\em IEEE Transactions on Automatic Control}, 62(12):6183--6194, 2017.

\bibitem{shoh21}
S.~Takata, Y.~Kato, and H.~Nakao.
\newblock Fast optimal entrainment of limit-cycle oscillators by strong periodic inputs via phase-amplitude reduction and {F}loquet theory.
\newblock {\em Chaos: An Interdisciplinary Journal of Nonlinear Science}, 31(9), 2021.

\bibitem{wigg03}
S.~Wiggins.
\newblock {\em Introduction to applied nonlinear dynamical systems and chaos}, volume~2.
\newblock Springer, 2003.

\bibitem{will15}
M.~O. Williams, I.~G. Kevrekidis, and C.~W. Rowley.
\newblock A data--driven approximation of the {k}oopman operator: Extending dynamic mode decomposition.
\newblock {\em Journal of Nonlinear Science}, 25(6):1307--1346, 2015.

\bibitem{wils20ddred}
D.~Wilson.
\newblock A data-driven phase and isostable reduced modeling framework for oscillatory dynamical systems.
\newblock {\em Chaos: An Interdisciplinary Journal of Nonlinear Science}, 30(1):013121, 2020.

\bibitem{wils20highacc}
D.~Wilson.
\newblock Phase-amplitude reduction far beyond the weakly perturbed paradigm.
\newblock {\em Physical Review E}, 101(2):022220, 2020.

\bibitem{wils21input}
D.~Wilson.
\newblock Analysis of input-induced oscillations using the isostable coordinate framework.
\newblock {\em Chaos: An Interdisciplinary Journal of Nonlinear Science}, 31(2):023131, 2021.

\bibitem{wils21dd}
D.~Wilson.
\newblock Data-driven inference of high-accuracy isostable-based dynamical models in response to external inputs.
\newblock {\em Chaos: An Interdisciplinary Journal of Nonlinear Science}, 31(6):063137, 2021.

\bibitem{wils21adapt}
D.~Wilson.
\newblock An adaptive phase-amplitude reduction framework without $\mathcal{O}(\epsilon)$ constraints on inputs.
\newblock {\em SIAM Journal on Applied Dynamical Systems}, 21(1):204--230, 2022.

\bibitem{wils25}
D.~Wilson.
\newblock Identification and computation of slow manifolds using the isostable coordinate system.
\newblock {\em arXiv preprint:2507.13997}, 2025.

\bibitem{wils17isored}
D.~Wilson and B.~Ermentrout.
\newblock Greater accuracy and broadened applicability of phase reduction using isostable coordinates.
\newblock {\em Journal of Mathematical Biology}, 76(1-2):37--66, 2018.

\bibitem{wils16isos}
D.~Wilson and J.~Moehlis.
\newblock Isostable reduction of periodic orbits.
\newblock {\em Physical Review E}, 94(5):052213, 2016.

\bibitem{winf01}
A.~Winfree.
\newblock {\em The Geometry of Biological Time}.
\newblock Springer Verlag, New York, second edition, 2001.

\end{thebibliography}
\end{document}